\documentclass[12inch,a4paper]{article}
\usepackage[top=3cm, bottom=2.5cm, left=2.4cm, right=2.7cm]{geometry}
\usepackage{graphicx}
\usepackage{epstopdf}
\usepackage{booktabs}
\usepackage{blkarray}
\usepackage[british]{babel}

 \usepackage[sort,comma,authoryear]{natbib} 

    \setcitestyle{citesep={;}}
\setcitestyle{square}

\usepackage{amssymb}
\usepackage{setspace}
\usepackage{amsfonts}
\usepackage{amsmath}
\usepackage{float}
\usepackage{color}
\usepackage[dvipsnames]{xcolor}
 \usepackage{colortbl}
\usepackage{url}

\definecolor{Crimson}{rgb}{0.6471, 0.1098, 0.1882}

\usepackage{epstopdf}
\newtheorem{lemma}{Lemma}
\newtheorem{remark}{\normalfont \textit{Remark}}

\newtheorem{proposition}{Proposition}
\newtheorem{corollary}{Corollary}

\usepackage{tcolorbox}
\usepackage{epstopdf}
\usepackage{subcaption}
\providecommand{\keywords}[1]{\textit{Key words:} #1}

\usepackage[labelfont=bf]{caption}

\usepackage{titlesec}
\titleformat*{\subsection}{\normalfont \large}
\usepackage{abstract}
\usepackage{slashbox}
\begin{document}



\title{Asymptotic properties of a componentwise ARH(1) plug-in predictor}

\author{Javier \'Alvarez-Li\'ebana$^1$, Denis Bosq$^2$ and M. Dolores Ruiz--Medina$^1$}
\maketitle
\begin{flushleft}
$^1$ Department of Statistics and O. R., University of Granada, Spain. 
$^2$ LSTA, Universit\'e Pierre et Marie Curie--Paris 6, Paris, France.

\textit{E-mail: javialvaliebana@ugr.es}
\end{flushleft}

\doublespacing



\renewcommand{\absnamepos}{flushleft}
\setlength{\absleftindent}{0pt}
\setlength{\absrightindent}{0pt}
\renewcommand{\abstractname}{Summary}
\begin{abstract}
This paper presents new results on prediction of linear processes in function spaces. The autoregressive Hilbertian process framework of order one (ARH(1) process framework) is adopted. A componentwise estimator of the autocorrelation operator is formulated, from the moment--based estimation of its diagonal coefficients, with respect to the orthogonal eigenvectors of the auto-covariance operator, which are assumed to be known. Mean-square convergence to the theoretical autocorrelation operator, in the space of Hilbert-Schmidt operators, is proved. Consistency then follows in that space. For the associated ARH(1) plug-in predictor, mean absolute convergence to the corresponding conditional expectation, in the considered Hilbert space, is obtained. Hence, consistency in that space also holds. A simulation study is undertaken to illustrate the finite-large sample behavior of the formulated componentwise estimator and predictor. The performance of the presented approach is compared with alternative approaches in the previous and current ARH(1) framework literature, including the case of unknown eigenvectors.

\vspace{0.5cm}
\textbf{Journal of Multivariate Analysis, 155, pp. 12-34.}
 DOI: \url{doi.org/10.1016/j.jmva.2016.11.009}
\end{abstract}


\textcolor{Crimson}{\section{Introduction}
\label{A3:sec:1}}

In the last few decades, an extensive literature on statistical inference from functional random variables has emerged.
This work was motivated in part by the statistical analysis of high--dimensional data, as well as data of a continuous (infinite-dimensional)
nature; see, e.g., \cite{Bosq00,Bosq07,DedeckerMerlevede03,FerratyVieu06,Merlevede96b,Merlevede97,RamsaySilverman05,Ruiz12}. New developments in functional data analysis are described,
e.g., in \cite{Bongiornoetal14,Cuevas14,HorvathKokoszka12,HsingEubank15}, and in a recent Special Issue of this journal \cite{GoiaVieu16}.

 These references include  a nice summary on the
statistics theory for functional data, contemplating  covariance operator theory and
eigenfunction expansion, perturbation theory, smoothing and regularization,  probability measures on a Hilbert spaces,  functional principal
component analysis,  functional counterparts of the  multivariate canonical correlation analysis, the two sample problem and the change point
problem, functional linear models, functional test for independence, functional time series theory, spatially distributed curves,
 software packages and   numerical implementation of the statistical procedures discussed,  among other topics.
 
The special case of functional regression models, in which the predictor is a random function and the response is
scalar, has been particularly well studied. Various specifications of the functional regression parameter arise in fields such as biology, climatology, chemometrics, and economics. To avoid the computational (high--dimensional) limitations of the
nonparametric approach, several parametric and semi--parametric methods have been proposed; see, e.g., \cite{Ferratyetal12} and the
references therein.  In \cite{Ferratyetal12}, a combination of a spline approximation and the one--dimensional Nadaraya--Watson approach
was proposed to avoid high dimensionality issues. Generalizations to the case of more regressors (all functional, or both
functional and real) were also addressed in the nonparametric, semi--parametric, and parametric frameworks; for an
overview, see \cite{AneirosVieu06,FebreroManteiga13,FerratyVieu09}.

In the nonparametric regression framework, the case where the covariates and the response are functional was
considered by \cite{Ferratyetal12}, where a functional version of the Nadaraya--Watson estimator was proposed for the
estimation of the regression operator and shown to be point--wise asymptotically normal. Resampling techniques were used
to overcome the difficulties arising in the estimation of the asymptotic bias and variance.
Semi--functional partial linear regression, introduced in \cite{AneirosVieu08}, allows the prediction of a real-valued random variable from
a set of real--valued explanatory variables, and a time--dependent functional explanatory variable. Motivated by genetic and
environmental applications, a semi--parametric maximum likelihood method for  the estimation of odds ratio association
parameters was developed by \cite{Chenetal12} in a high--dimensional data context.

In the autoregressive Hilbertian time series framework, several estimation and prediction procedures have been
proposed and studied. \cite{Mas99} established, under suitable conditions, the asymptotic normal distribution of the formulated
estimator of the autocorrelation operator, based on projection into the theoretical eigenvectors. In \cite{Bosq00,BosqBlanke07}, the problem of
prediction of linear processes in function spaces was addressed. In particular, sufficient conditions for the consistency of
the empirical autocovariance and cross--covariance operators were obtained. The asymptotic normal distribution of the
empirical autocovariance operator was also derived. Moreover, the asymptotic properties of the empirical eigenvalues and
eigenvectors were analysed.

\cite{Guillas01} established the efficiency of a componentwise estimator of the autocorrelation operator, based on projection
into the empirical eigenvector system of the autocovariance operator. Consistency, in the space of bounded linear operators,
of the formulated estimator of the autocorrelation operator, and of its associated ARH(1) plug--in predictor was later proved
by \cite{Mas04}. He derived sufficient conditions for the weak convergence of the ARH(1) plug--in predictor to a Hilbert--valued
Gaussian random variable (see \cite{Mas07}). Simultaneously, \cite{MasMennetau03a} obtained high deflection results or large and moderate deviations
for infinite--dimensional autoregressive processes. Furthermore, the law of the iterated logarithm for the covariance
operator estimator was formulated by \cite{Menneteau05}.

The main properties of the class of autoregressive Hilbertian processes with random coefficients were investigated by
\cite{Mourid04}. \cite{Kargin08} gave interesting extensions of the autoregressive Hilbertian framework, based on the
spectral decomposition of the autocorrelation operator, and not of the autocovariance operator. The first generalization
on autoregressive processes of order greater than one was proposed by \cite{Mourid93}, in order to improve prediction.
ARHX(1) models; i.e., autoregressive Hilbertian processes with exogenous variables were studied by \cite{DamonGuillas02,DamonGuillas05}.
In \cite{Guillas00,Guillas01} a doubly stochastic formulation of the autoregressive Hilbertian process was investigated. The ARHD model
was introduced by \cite{MarionPumo04}, taking into account the regularity of trajectories through the derivatives. The
conditional autoregressive Hilbertian process (CARH process) was considered by \cite{Cugliari11}, developing parallel projection
estimation methods to predict such processes. In the Banach--valued context, we refer to the papers by \cite{BensmainMourid01,DehlingSharipov05,Pumo92,Pumo98}, among
others.

In this paper, we assume that the autocorrelation operator belongs to the Hilbert--Schmidt class, and admits a diagonal
spectral decomposition in terms of the orthogonal eigenvector system of the autocovariance operator. Such is the case, e.g., of
an autocorrelation operator defined as a continuous function of the autocovariance operator. A componentwise estimator
of the autocorrelation operator is then constructed in terms of the eigenvectors of the autocovariance operator, which are
assumed to be known. This occurs when the random initial condition is defined as the solution, in the mean--square sense, of a
stochastic differential equation driven by white noise. Beyond this case, the sparse representation and whitening properties
of wavelet bases can be exploited to obtain a diagonal representation of the autocovariance and cross--covariance operators,
in terms of a common and known wavelet basis. Unconditional bases, like wavelet bases, also allow the diagonal spectral
series representation of the distributional kernels of Calder\'on-Zygmund operators.

Under the assumptions stated in \textcolor{Crimson}{Appendices} \ref{A3:sec:2}--\ref{A3:gc}, we establish the convergence in the $\mathcal{L}^{2}$-sense of a componentwise
estimator of the autocorrelation operator in the space of Hilbert--Schmidt operators $\mathcal{S} \left( H \right),$ i.e., $\mathcal{L}^{2}_{\mathcal{S} \left( H \right)} \left(\Omega,
\mathcal{A}, \mathcal{P} \right),$ is derived. Consistency then follows in $\mathcal{S} \left( H \right)$. Under the same conditions, consistency in H of the associated ARH(1) plug--in predictor is obtained, from its convergence in the $\mathcal{L}^{1}$-sense in the Hilbert space $H,$ i.e., in the space  $\mathcal{L}^{1}_{H } \left(\Omega, \mathcal{A}, \mathcal{P} \right)$. The Gaussian framework is analysed in \textcolor{Crimson}{Appendix} \ref{A3:gc} and illustrated in \textcolor{Crimson}{Appendix} \ref{A3:sec:4}, where examples show the behaviour of
the proposed componentwise autocorrelation operator estimator, and associated predictor, for large sample sizes. We also
present there a comparative study with alternative ARH(1) prediction techniques, including componentwise parameter
estimation of the autocorrelation operator, from known and unknown eigenvectors, as well as kernel (nonparametric)
functional estimation, and penalized, spline and wavelet, estimation. Final comments on the application of the proposed
approach from real data are provided in \textcolor{Crimson}{Appendix} \ref{A3:sec:5}.

%
%

\textcolor{Crimson}{\section{Preliminaries}
\label{A3:sec:2}}

This section contains the preliminary definitions and lemmas that will be used to derive the main results of this paper. In the following, $H$ denotes a real separable Hilbert space. Recall that, from \cite{Bosq00},  a zero--mean
ARH(1) process $X= \left\lbrace X_{n},\ n \in \mathbb{Z} \right\rbrace$ satisfies, for all $n \in \mathbb{Z}$, the equation

\begin{equation}
X_{n } = \rho \left(X_{n-1} \right) + \varepsilon_{n},\label{A3:24bb}
\end{equation}

\noindent  where $\rho $ denotes the autocorrelation operator of the
process $X,$ which belongs to the space $\mathcal{L}(H)$ of bounded
linear operators, such that $\Vert \rho^{k}
\Vert_{\mathcal{L}\left(H\right)}<1,$  for all integers $k\geq k_{0}$ beyond a certain $k_{0} \geq 1$, with $\|\cdot\|_{\mathcal{L}(H)}$ denoting the
norm in the space $\mathcal{L}(H).$ The Hilbert--valued innovation
process
 $\varepsilon= \left\lbrace \varepsilon_{n }, \ n\in
\mathbb{Z} \right\rbrace$ is assumed to be a strong--white noise which is uncorrelated with
the random initial condition. That is, $\varepsilon$ is a Hilbert--valued zero--mean stationary process, with independent and identically
distributed components in time, with
$\sigma^{2}_{\varepsilon}={\rm E} \left\lbrace \|\varepsilon_{n}\|_{H}^{2} \right\rbrace <\infty,$ for
all $n\in \mathbb{Z}.$ We restrict our attention here  to the case
where $\rho$ is such that  $$\|\rho\|_{\mathcal{L}(H)}<1.$$

The following assumptions are made.

\bigskip

\noindent \textcolor{Aquamarine}{\textbf{Assumption A1.}} The autocovariance operator $$C= {\rm E} \left\lbrace X_{n}\otimes
X_{n} \right\rbrace = {\rm E} \left\lbrace X_{0}\otimes X_{0} \right\rbrace, \quad n \in \mathbb{Z},$$ is
a positive, self--adjoint and trace operator. As a result, it admits  the following diagonal spectral representation
\begin{equation}
C=\displaystyle \sum_{j=1}^{\infty}C_{j}\phi_{j}\otimes \phi_{j}, \nonumber 
\end{equation}

\noindent in terms of an orthonormal system $\left\lbrace \phi_{j}, \ j\geq 1 \right\rbrace$ of eigenvectors which are known. Here, $$C_1 \geq C_2 \geq \dots \geq
C_j \geq  \dots
> 0$$ denote the real
positive eigenvalues of $C$ arranged in decreasing order of magnitude and
$$\displaystyle
\sum_{j=1}^{\infty} C_j < \infty.$$

\bigskip

\noindent \textcolor{Aquamarine}{\textbf{Assumption A2.}} The autocorrelation operator $\rho$ is a self--adjoint and Hilbert--Schmidt operator, admitting the diagonal
spectral decomposition
\begin{equation}
\rho = \displaystyle \sum_{j=1}^{\infty} \rho_j \phi_j \otimes
\phi_j,\quad \displaystyle \sum_{j=1}^{\infty} \rho_{j}^{2} <
\infty, \nonumber 
\end{equation}

\noindent where $\left\lbrace \rho_j, \ j \geq 1 \right\rbrace$ is the system of
eigenvalues of the autocorrelation operator $\rho,$ with respect to
the orthonormal system  of eigenvectors $\left\lbrace \phi_j, \ j \geq 1
\right \rbrace$ of the autocovariance operator $C$.

\bigskip

Note that, under \textcolor{Aquamarine}{\textbf{Assumption A2}}, 
$$\|\rho\|_{\mathcal{L}(H)}=\displaystyle \sup_{j\geq 1}\left| \rho_{j} \right| <1.$$

\bigskip

\begin{remark}
\label{A3:rem1def}
\textit{\textcolor{Aquamarine}{\textbf{Assumption A2}} holds, in particular, when operator $\rho $  is defined as a continuous function of operator C (see \cite[pp. 119--140]{DautrayLions90} and \textcolor{Crimson}{Remark} \ref{A3:remark3def}).}
\end{remark}

\bigskip

In the following, for any $n\in \mathbb{Z},$ let  $$D= {\rm E} \left\lbrace X_{n}\otimes
X_{n+1} \right\rbrace = {\rm E} \left\lbrace X_{0}\otimes X_{1} \right\rbrace$$  be the cross--covariance operator of
the ARH(1) process $X$.

\bigskip

\begin{remark}
\textit{Under \textcolor{Aquamarine}{\textbf{Assumptions A1--A2}}, it follows from equation (\ref{A3:24bb}) that
 \begin{equation}
 C_{\varepsilon} = C_ \rho C \rho = \displaystyle \sum_{j=1}^{\infty} C_j \left(1 - \rho_{j}^{2} \right) \phi_j \otimes \phi_j = \displaystyle \sum_{j=1}^{\infty} \sigma_{j}^{2} \phi_j \otimes \phi_j . \nonumber
 \end{equation}}
 \end{remark}

\bigskip

By projecting equation (\ref{A3:24bb})  into the orthonormal system $\left\lbrace \phi_j, \ j \geq 1 \right\rbrace$, we also have, for each $ j \geq 1$ and all $n \in \mathbb{Z}$, the AR(1) equation
\begin{equation}
X_{n,j}=\rho_{j}X_{n-1,j}+\varepsilon_{n,j},\quad n\in \mathbb{Z},
\label{A3:14}
\end{equation}

\noindent where $X_{n,j}=\left\langle
X_{n},\phi_{j}\right\rangle_{H}$ and $\varepsilon_{n,j}=\left\langle
\varepsilon_{n},\phi_{j}\right\rangle_{H},$ for all $n\in
\mathbb{Z}$. From equation (\ref{A3:14}), we have, for each $j\geq 1$ and all $n \in \mathbb{Z}$,
\begin{eqnarray}
\rho_{j}  &=& \rho (\phi_{j})(\phi_{j})=\left\langle \phi_{j},DC^{-1}(\phi_{j})\right\rangle_{H} = \left\langle D(\phi_{j}),\phi_{j}\right\rangle_{H}\left\langle C^{-1}(\phi_{j}),\phi_{j}\right\rangle_{H}\nonumber\\
&=&\frac{{\rm E} \left\lbrace X_{n,j}X_{n-1,j}\right\rbrace}{{\rm E} \left\lbrace X_{n-1,j}^{2}
\right\rbrace} = \frac{D_j}{C_j}, \quad  n\in \mathbb{Z},
 \label{A3:13}
\end{eqnarray}
\noindent where   $$D_j=\left\langle D(\phi_{j}),\phi_{j}\right\rangle_{H}= {\rm E} \left\lbrace X_{n,j}X_{n-1,j}\right\rbrace, \quad C_j^{-1}=[{\rm E} \left\lbrace X_{n-1,j}^{2}\right\rbrace]^{-1}, \quad X_{n,j}=\left\langle X_{n},\phi_{j}\right\rangle_{H},$$ given that, for all $j \geq 1$,
\begin{equation}
D=\sum_{j=1}^{\infty }D_{j}\phi_{j}\otimes \phi_{j},\quad
D_{j}=\rho_{j}C_{j},\quad j\geq 1.\label{A3:eqproyar1}
\end{equation}

Let us now  consider the  Banach space $L_{\mathcal{H}}^{2} \left( \Omega,
\mathcal{A}, \mathcal{P} \right)$ of the  equivalence classes of $\mathcal{L}_{\mathcal{H}}^{2} \left( \Omega,
\mathcal{A}, \mathcal{P} \right),$ the space of zero--mean
second--order Hilbert--valued random variables ($\mathcal{H}$--valued
random variables) with finite seminorm given by
\begin{equation}
\left\| Z \right\|_{\mathcal{L}_{\mathcal{H}}^{2} \left( \Omega,
\mathcal{A}, \mathcal{P} \right)}=\sqrt{{\rm E} \left\lbrace \left\| Z
\right\|_{\mathcal{H}}^{2} \right\rbrace} ,\quad \forall Z\in
\mathcal{L}_{\mathcal{H}}^{2} \left( \Omega, \mathcal{A},
\mathcal{P} \right). \nonumber 
\end{equation}

\noindent That is, for $Z,Y\in \mathcal{L}_{\mathcal{H}}^{2} \left( \Omega,
\mathcal{A}, \mathcal{P} \right),$  $Z$ and $Y$ belong to the same equivalence class if and only if  $${\rm E} \left\lbrace \left\| Z-Y
\right\|_{\mathcal{H}} \right\rbrace =0.$$

The convergence in the seminorm of $\mathcal{L}_{\mathcal{S}(H)}^{2} \left( \Omega, \mathcal{A},\mathcal{P} \right)$ will be considered in \textcolor{Crimson}{Proposition}  \ref{A3:proposition1}, where $\mathcal{H}=\mathcal{S}(H)$ denotes the Hilbert
space of Hilbert--Schmidt operators on a Hilbert space $H$.

For each $n\in \mathbb{Z},$ let us consider  the following
biorthogonal representation of the functional value $X_{n}$ of the
ARH(1) process $X = \left\lbrace X_n, \ n \in \mathbb{Z} \right\rbrace$, and of the functional value $\varepsilon_{n}$ of its innovation process:

\begin{eqnarray}
X_{n}&=&\sum_{j=1}^{\infty }\sqrt{C_{j}}\frac{\left\langle
X_{n},\phi_{j}\right\rangle_{H}}{\sqrt{C_{j}}}\phi_{j}=\sum_{j=1}^{\infty
}\sqrt{C_{j}}\eta_{j}(n)\phi_{j}, \label{A3:15}\\
\varepsilon_{n}&=&\sum_{j=1}^{\infty }\sigma_{j}\frac{\left\langle
\varepsilon_{n},\phi_{j}\right\rangle_{H}}{\sigma_{j}}\phi_{j}=\sum_{j=1}^{\infty
}\sigma_{j}\widetilde{\eta}_{j}(n)\phi_{j}, \label{A3:16}
\end{eqnarray}

\noindent where
$$\eta_j (n) = \frac{\left\langle
X_{n},\phi_{j}\right\rangle_{H}}{\sqrt{C_{j}}} =
\frac{X_{n,j}}{\sqrt{C_{j}}}, \quad \widetilde{\eta}_{j} (n) =
\frac{\left\langle
\varepsilon_{n},\phi_{j}\right\rangle_{H}}{\sigma_{j}} =
\frac{\varepsilon_{n,j}}{\sigma_{j}}, \quad n \in
\mathbb{Z},~j \geq 1.$$ Here, under  \textcolor{Aquamarine}{\textbf{Assumptions A1--A2}}, for $C_{\varepsilon}= {\rm E} \left\lbrace \varepsilon_{n}\otimes \varepsilon_{n} \right\rbrace= {\rm E} \left\lbrace\varepsilon_{0}\otimes \varepsilon_{0} \right\rbrace,~n\in \mathbb{Z},$
$$C_{\varepsilon} \left( \phi_{j} \right) =\sigma_{j}^{2}\phi_{j},\quad j\geq 1,$$
where, as before, $\left\lbrace \phi_{j},\ j \geq 1 \right\rbrace$ denotes the system of eigenvectors of the autocovariance operator $C,$ and  $$\displaystyle \sum_{j=1}^{\infty} \sigma_{j}^{2}=\sigma^{2}_{\varepsilon}= {\rm E} \left\lbrace \|\varepsilon_{n}\|_{H}^{2} \right\rbrace,$$ for all $n\in \mathbb{Z}.$

The following lemma provides the convergence, in the seminorm of $\mathcal{L}_{H}^{2}(\Omega,\mathcal{A},\mathcal{P}),$ of the
 series expansions (\ref{A3:15})--(\ref{A3:16}).

\bigskip

\begin{lemma}
\label{A3:proposition1} Let $X =  \left\lbrace X_n, \ n \in \mathbb{Z} \right\rbrace$ be a
zero--mean ARH(1) process. Under \textcolor{Aquamarine}{\textbf{Assumptions A1--A2}}, for
any $n\in \mathbb{Z},$ the following limit holds
\begin{equation}
\displaystyle \lim_{M \to \infty}{\rm E} \left\lbrace \left\| X_n -
\widehat{X}_{n,M}\right\|_{H}^{2} \right\rbrace = 0, \nonumber 
\end{equation}

\noindent where $\widehat{X}_{n,M} = \displaystyle \sum_{j=1}^{M} \sqrt{C_j}
\eta_j (n) \phi_j$. Furthermore,
\begin{equation}
\displaystyle \lim_{M \to \infty} \left\| {\rm E} \left\lbrace \left(X_n -
\widehat{X}_{n,M}\right) \otimes \left(X_n -
\widehat{X}_{n,M}\right)\right\rbrace \right\|_{\mathcal{S} \left( H
\right)}^{2} = 0. \nonumber 
\end{equation}

Similar assertions hold for the biorthogonal series representation
$$\varepsilon_{n}=\sum_{j=1}^{\infty }\sigma_{j}\frac{\left\langle
\varepsilon_{n},\phi_{j}\right\rangle_{H}}{\sigma_{j}}\phi_{j}=\sum_{j=1}^{\infty
}\sigma_{j}\widetilde{\eta}_{j}(n)\phi_{j}.$$
\end{lemma}

\begin{proof}

Under \textcolor{Aquamarine}{\textbf{Assumption A1}}, from the trace property of $C,$ the
sequence $$\left\lbrace \widehat{X}_{n,M} = \displaystyle \sum_{j=1}^{M}
\sqrt{C_j} \eta_j (n) \phi_j,~M\geq 1\right\rbrace$$ satisfies, for $M$  sufficiently large, and $L>0,$ arbitrary,
\begin{eqnarray} \|\widehat{X}_{n,M+L}-\widehat{X}_{n,M}\|^{2}_{\mathcal{L}^{2}_{H}(\Omega, \mathcal{A},
P)} &=& {\rm E} \left\lbrace \|\widehat{X}_{n,M+L}-\widehat{X}_{n,M}\|^{2}_{H} \right\rbrace \nonumber\\
& =& \sum_{j=M+1}^{M+L}\sum_{k=M+1}^{M+L}\sqrt{C_j}\sqrt{C_k}
{\rm E} \left\lbrace \eta_j (n)\eta_{k}(n) \right\rbrace \left\langle\phi_j,\phi_{k}\right\rangle_{H} \nonumber \\
&=&\sum_{j=M+1}^{M+L}C_j \rightarrow 0, \quad \text{when } M\rightarrow \infty, \label{A3:inq}
\end{eqnarray}

\noindent since, under \textcolor{Aquamarine}{\textbf{Assumption A1}},  $\displaystyle \sum_{j=1}^{\infty}C_{j}<\infty$. Hence, $\left\lbrace \displaystyle \sum_{j=1}^{M}C_{j}, \ M\geq 1 \right\rbrace$ is a Cauchy sequence. Thus, $$\displaystyle \lim_{M \to \infty} \displaystyle \sum_{j=M+1}^{M+L}C_j = 0,$$ 
for $L>0$ arbitrary. From equation (\ref{A3:inq}), $$\left\lbrace\widehat{X}_{n,M} = \displaystyle \sum_{j=1}^{M}
\sqrt{C_j} \eta_j (n) \phi_j,\ M\geq 1\right\rbrace$$ is also a Cauchy sequence in $\mathcal{L}_{H}^{2}(\Omega,\mathcal{A},P).$ Thus, the sequence $\left\lbrace \widehat{X}_{n,M}, \ M\geq 1 \right\rbrace$ has
finite limit in \linebreak $\mathcal{L}_{H}^{2}(\Omega, \mathcal{A}, \mathcal{P})$, for all $n \in \mathbb{Z}$.

  Furthermore,
\begin{eqnarray}
\displaystyle \lim_{M \to \infty} {\rm E} \left\lbrace \left\| X_n -
\widehat{X}_{n,M}\right\|_{H}^{2} \right\rbrace & =  & {\rm E} \left\lbrace \left\| X_n
\right\|_{H}^{2} \right\rbrace  +  \displaystyle \lim_{M \to \infty}   \displaystyle \sum_{j=1}^{M}  \displaystyle \sum_{h=1}^{M} \sqrt{C_j} \sqrt{C_h} {\rm E} \left\lbrace  \eta_j (n)\eta_h (n) \right\rbrace \langle \phi_j,\phi_h \rangle_H  \nonumber \\
& - & 2 \displaystyle \lim_{M \to \infty} \displaystyle
\sum_{j=1}^{M} \sqrt{C_j} {\rm E} \left\lbrace\langle X_n,\eta_j (n) \phi_j
\rangle_H \right\rbrace = \sigma_{X}^{2} \nonumber \\
&-& \displaystyle \lim_{M \to \infty} \displaystyle \sum_{j=1}^{M} C_j =0 . \nonumber\\
\label{A3:19}
\end{eqnarray}

In the derivation of the identities in (\ref{A3:inq})--(\ref{A3:19}),
we have applied that, for every $j,~h\geq 1,$

\begin{eqnarray}
C \left(\phi_{j} \right) &=&C_{j}\phi_{j}, \quad \quad \sigma_{X}^{2}= {\rm E} \left\lbrace \|X_{n}\|^{2}_{H} \right\rbrace = \sum_{j=1}^{\infty}C_{j}<+\infty,\quad \left\langle \phi_{j},\phi_{h}\right\rangle_{H} = \delta_{j,h}, \nonumber\\
{\rm E} \left\lbrace \eta_{j}(n)\eta_{h}(n)\right\rbrace &=&\delta_{j,h},\quad 
{\rm E} \left\lbrace \left\langle X_{n}, \eta_{j}(n)\phi_{j}\right\rangle_{H}\right\rbrace =  \sqrt{C_{j}}. \nonumber \\\label{A3:20}
\end{eqnarray}

Moreover, from identities in (\ref{A3:20}),
\begin{eqnarray}
&&\left\| {\rm E} \left\lbrace \left(X_{n}-\lim_{M\rightarrow \infty}
\widehat{X}_{n,M} \right)\otimes \left(X_{n}-\lim_{M\rightarrow
\infty} \widehat{X}_{n,M} \right) \right\rbrace
\right\|^{2}_{\mathcal{S}(H)}
\nonumber\\
&=&\left\| {\rm E} \left\lbrace X_{n}\otimes X_{n} \right\rbrace + \displaystyle \lim_{M\to \infty}
\displaystyle \sum_{j=1}^{M}
\displaystyle\sum_{h=1}^{M}\sqrt{C_{j}}\sqrt{C_{h}}\phi_{j}\otimes\phi_{h}{\rm E} \left\lbrace\eta_{j}(n)\eta_{h}(n) \right\rbrace \right.\nonumber\\
& &\left.- 2\lim_{M\rightarrow \infty}\sum_{j=1}^{M} {\rm E} \left\lbrace X_{n}\otimes
\sqrt{C_{j}}\eta_{j}(n)\phi_{j} \right\rbrace \right\|^{2}_{\mathcal{S}(H)}\nonumber \\
&=& \left\| {\rm E} \left\lbrace X_{n}\otimes X_{n} \right\rbrace+\lim_{M\rightarrow
\infty}\left[\sum_{j=1}^{M}C_{j}\phi_{j}\otimes
\phi_{j}-2\sum_{j=1}^{M}C_{j}\phi_{j}\otimes \phi_{j}\right]\right\|^{2}_{\mathcal{S}(H)}\nonumber\\
&=&\left\| {\rm E} \left\lbrace X_{n}\otimes X_{n} \right\rbrace-\lim_{M\rightarrow
\infty}\sum_{j=1}^{M}C_{j}\phi_{j}\otimes
\phi_{j}\right\|^{2}_{\mathcal{S}(H)}=0.\label{A3:21}
\end{eqnarray}

In a similar way, we can derive the convergence to $\varepsilon_{n},$ in
$\mathcal{L}_{H}^{2}(\Omega,\mathcal{A},\mathcal{P}),$ of the series
$\displaystyle \sum_{j=1}^{\infty
}\sigma_{j}\widetilde{\eta}_{j}(n)\phi_{j},$  for every $n\in
\mathbb{Z},$ since $\varepsilon $ is assumed to be  strong--white
noise, and hence,  its covariance operator $C_{\varepsilon}$  is
 in the trace class.   We
can also obtain an analogous to equation (\ref{A3:21}).

 \hfill \hfill \textcolor{Aquamarine}{$\blacksquare$}
\end{proof}

\bigskip

In equations (\ref{A3:15})--(\ref{A3:16}),  for every $n\in \mathbb{Z},$
\begin{eqnarray}
& & {\rm E} \left\lbrace \eta_{j}(n) \right\rbrace =0,\quad {\rm E} \left\lbrace \eta_{j}(n)\eta_{h}(n) \right\rbrace =\delta_{j,h},\quad j,h\geq 1, \quad n\in \mathbb{Z},\label{A3:22}\\
& & {\rm E} \left\lbrace\widetilde{\eta}_{j}(n) \right\rbrace=0,\quad
{\rm E} \left\lbrace \widetilde{\eta}_{j}(n)\widetilde{\eta}_{h}(n) \right\rbrace=\delta_{j,h},\quad
j,h\geq 1,\quad n\in \mathbb{Z}. \nonumber 
\end{eqnarray}

Note that, from \textcolor{Aquamarine}{\textbf{Assumption A2}} for each $j\geq 1,$ $ \left\lbrace X_{n,j}, \ n\in
\mathbb{Z} \right\rbrace$ in equation (\ref{A3:14}) defines a stationary and
invertible AR(1) process. In addition, from equations (\ref{A3:15}) and
(\ref{A3:20}),
 for every $n \in \mathbb{Z},$ \linebreak and  $j,p \geq 1,$
\begin{eqnarray}
X_{n} &=&  \displaystyle \sum_{j=1}^{\infty} X_{n,j} \phi_j, \nonumber\\
{\rm E} \left\lbrace X_{n,j} X_{n,p} \right\rbrace &=& \displaystyle
\sum_{k=0}^{\infty} \displaystyle \sum_{h=0}^{\infty}  \rho_{j}^{k}
\rho_{p}^{h} {\rm E} \left\lbrace \varepsilon_{n-k,j} \varepsilon_{n-h,p}
\right\rbrace = \delta_{j,p}
 \displaystyle \sum_{k=0}^{\infty} \rho_{j}^{2k} \sigma_{j}^{2} = \delta_{j,p} \frac{\sigma_{j}^{2}}{1 - \rho_{j}^{2}}, \nonumber\\
{\rm E} \left\lbrace \left\| X_{n} \right\|_{H}^{2} \right\rbrace &=& \displaystyle \sum_{j=1}^{\infty} {\rm E} \left\lbrace X_{n,j}^2 \right\rbrace = \displaystyle \sum_{j=1}^{\infty} \langle C \left( \phi_j \right), \phi_j \rangle_H = \displaystyle \sum_{j=1}^{\infty} C_j = \sigma_{X}^{2} < \infty,  \nonumber \\ \label{A3:25a}
\end{eqnarray}

\noindent which implies that $$C_j = \frac{\sigma_{j}^{2}}{1 -
\rho_{j}^{2}}, \quad j\geq 1.$$   In particular, we obtain, for
each $j \geq 1,$ and for every $n \in \mathbb{Z},$
\begin{eqnarray}
{\rm E} \left\lbrace \eta_j (n) \eta_j (n+1) \right\rbrace &=& {\rm E} \left\lbrace \frac{X_{n,j}}{\sqrt{C_j}} \frac{X_{n+1,j}}{\sqrt{C_j}} \right\rbrace = \frac{{\rm E} \left\lbrace X_{n,j} X_{n+1,j} \right\rbrace}{C_j}
\nonumber \\
 &=& \frac{\displaystyle \sum_{k=0}^{\infty} \displaystyle \sum_{h=0}^{\infty}  \rho_{j}^{k+h} {\rm E} \left\lbrace \varepsilon_{n-k,j}\varepsilon_{n+1-h,j} \right\rbrace}{C_j}  \nonumber \\
 &=& \frac{\displaystyle \sum_{k=0}^{\infty}   \rho_{j}^{2k + 1}
\sigma_{j}^{2}}{C_j} = \frac{ \sigma_{j}^{2}}{C_j}
\frac{\rho_j}{1-\rho_{j}^{2}}=\rho_j.\label{A3:24a}
\end{eqnarray}

\bigskip

\begin{remark}
\label{A3:remark1}
\textit{From equation (\ref{A3:14}) and \textcolor{Crimson}{Lemma} \ref{A3:proposition1}, keeping in
mind that $$C_j = \frac{\sigma_{j}^{2}}{1 - \rho_{j}^{2}}, \quad j\geq 1,$$ the  following invertible and stationary AR(1) process
can be defined:
\begin{equation}
\eta_{j}(n)=\rho_{j}\eta_{j}(n-1)+\sqrt{1-\rho_{j}^{2}}\widetilde{\eta}_{j}(n),\quad
0< \rho_{j}^{2}\leq \rho_{j} <1, \label{A3:AE1}
\end{equation}
\noindent where, for each $j\geq 1,$ $\left\lbrace \eta_{j}(n),\ n\in
\mathbb{Z} \right\rbrace$ and $ \left\lbrace \widetilde{\eta}_{j}(n),\ n\in \mathbb{Z} \right\rbrace$ are
respectively  introduced in equations (\ref{A3:15})-(\ref{A3:16}). In
the following,  for each $j\geq 1,$ we assume that
$${\rm E} \left\lbrace \left(\widetilde{\eta}_{j}(n)\right)^{4}\right\rbrace <\infty, \quad n\in \mathbb{Z},$$ \noindent to ensure ergodicity for all
second--order moments, in the mean--square sense; see, e.g.,
\cite[pp. 192--193]{Hamilton94}.}
\end{remark}

\bigskip

Furthermore,
\begin{eqnarray}
D={\rm E}\left\lbrace X_{n}\otimes X_{n+1} \right\rbrace &=&\sum_{j=1}^{\infty }\sum_{p=1}^{\infty }{\rm E}\left\lbrace \left\langle
X_{n},\phi_{j}\right\rangle_{H}\left\langle
X_{n+1},\phi_{p}\right\rangle_{H}\right\rbrace \phi_{j}\otimes \phi_{p}\nonumber\\
&=&\sum_{j=1}^{\infty }\sum_{p=1}^{\infty }\sqrt{C_{j}}\sqrt{C_{p}}\frac{{\rm E}\left\lbrace \left\langle
X_{n},\phi_{j}\right\rangle_{H}\left\langle
X_{n+1},\phi_{p}\right\rangle_{H}\right\rbrace}{\sqrt{C_{j}}\sqrt{C_{p}}}\phi_{j}\otimes \phi_{p}\nonumber\\
&=&\sum_{j=1}^{\infty }\sum_{p=1}^{\infty }\sqrt{C_{j}}\sqrt{C_{p}}{\rm E}\left\lbrace \eta_{j}(n)\eta_{p}(n+1)\right\rbrace\phi_{j}\otimes \phi_{p}. \nonumber 
\end{eqnarray}

\bigskip

\begin{remark}
\label{A3:remark3def}
\textit{In particular,  \textcolor{Aquamarine}{\textbf{Assumption A2}} holds if  the following orthogonality condition is satisfied, for all $n \in \mathbb{Z}$ and $j,p \geq 1$,
\begin{eqnarray}
{\rm E}\left\lbrace \eta_{j}(n)\eta_{p}(n+1)\right\rbrace&=&\delta_{j,p},\nonumber 
\end{eqnarray}
\noindent where $\delta_{j,p}$ denotes the Kronecker Delta function. In practice, unconditional bases, e.g., wavelet bases, lead   to a sparse representation for functional data; see, e.g., \cite{Nason08,Ogden97,Vidakovic98} for statistically-oriented treatments.  Wavelet bases are also designed for sparse representation of kernels defining integral operators, in $L^{2}$ spaces with respect to a suitable measure (see \cite{Mallat09}). The Discrete Wavelet Transform (DWT) approximately decorrelates or \emph{whitens} data (see \cite{Vidakovic98}). In particular, operators   $C$ and $D$ could admit an almost  diagonal representation with respect to the self-tensorial tensorial product of a suitable  wavelet basis.}
\end{remark}

%
%

\textcolor{Crimson}{\section{Estimation and prediction results}
\label{A3:sec:3}}

 A componentwise estimator of the autocorrelation
operator and of the associated ARH(1)  plug--in predictor are
formulated in this section. Their convergence to the corresponding
theoretical functional values are derived in the
 spaces $\mathcal{L}^{2}_{\mathcal{S}(H)}(\Omega,\mathcal{A},\mathcal{P})$ and
$\mathcal{L}_{H}(\Omega,\mathcal{A},\mathcal{P}),$ respectively. Their consistency in the spaces $\mathcal{S}(H)$ and $H$
 then follows.

From equation (\ref{A3:13}), for each $j\geq 1,$ and for a given sample
size $n$, one can consider the
 usual  respective moment--based estimators $\widehat{D}_{n,j}$ and $\widehat{C}_{n,j}$ of $D_{j}$ and
 $C_{j},$ in the AR(1) framework, given by

\begin{eqnarray}
\widehat{D}_{n,j} &=& \frac{1}{n-1}\displaystyle \sum_{i=0}^{n-2}
X_{i,j}X_{i+1,j}, \quad \widehat{C}_{n,j} = \frac{1}{n}\displaystyle \sum_{i=0}^{n-1}
X_{i,j}^2. \nonumber 
\end{eqnarray}

The following truncated componentwise  estimator of $\rho$ is then
formulated:
\begin{equation}
\widehat{\rho}_{k_{n}}= \displaystyle \sum_{j=1}^{k_{n}}\widehat{\rho}_{n,j}\phi_{j}\otimes
\phi_{j}, \label{A3:24}
\end{equation}

\noindent where, for each $j\geq 1,$
\begin{equation}
\widehat{\rho}_{n,j} = \frac{\widehat{D}_{n,j}}{\widehat{C}_{n,j}} =
\frac{\frac{1}{n-1}\displaystyle \sum_{i=0}^{n-2}
X_{i,j}X_{i+1,j}}{\frac{1}{n}\displaystyle \sum_{i=0}^{n-1}
X_{i,j}^2} = \frac{n}{n-1}\frac{\displaystyle \sum_{i=0}^{n-2}
X_{i,j}X_{i+1,j}}{\displaystyle \sum_{i=0}^{n-1} X_{i,j}^2} .
\label{A3:25}
\end{equation}

 Here, the truncation parameter indicates that we
have considered the first $k_{n}$ eigenvectors associated with the
first $k_{n}$ eigenvalues, arranged in decreasing order of their
modulus magnitude. Furthermore, $k_{n}$ is such that
\begin{equation}
\displaystyle \lim_{n\rightarrow \infty }k_{n}=\infty,\quad \frac{k_{n}}{n}<
1,\quad n \geq 2.\label{A3:26}
\end{equation}

The following additional condition will be assumed on $k_{n}$ for
the derivation of the subsequent results:

\bigskip

\noindent \textcolor{Aquamarine}{\textbf{Assumption A3.}} The truncation parameter $k_{n}$
in (\ref{A3:24}) is such that $$\displaystyle \lim_{n \to
\infty}C_{k_n} \sqrt{n} = \infty.$$

\bigskip

\begin{remark}
\label{A3:remark2}
 \textcolor{Aquamarine}{\textbf{Assumption A3}} has also been considered  in \cite[p. 217]{Bosq00}, to ensure weak consistency of the proposed estimator of $\rho,$ as well as, in  \cite[Proposition 4]{Mas99}, in the derivation of asymptotic normality.
\end{remark}

\bigskip

From \textcolor{Crimson}{Remark} \ref{A3:remark1}, for each $j\geq 1,$
$\eta_{j}= \left\lbrace \eta_{j}(n),\ n\in \mathbb{Z} \right\rbrace$ in equation (\ref{A3:AE1})
defines a stationary and invertible AR(1) process, ergodic in the
mean--square sense; see, e.g., \cite{Bartlett46}. Therefore, in view of equations (\ref{A3:22}) and
(\ref{A3:24a}), for each $j \geq 1$, there exist two positive constants $K_{j,1}$ and
$K_{j,2}$ such that the following identities hold:
\begin{eqnarray}
&& \displaystyle \lim_{n \to \infty} \frac{{\rm E} \left\lbrace \left[1- \frac{1}{n} \displaystyle \sum_{i=0}^{n-1} \eta_{j}^{2} (i)\right]^{2}\right\rbrace }{\frac{1}{n}} = K_{j,1}, \label{A3:27} \\
&& \displaystyle \lim_{n \to \infty} \frac{{\rm E} \left\lbrace \left[ \rho_{j} -
\frac{1}{n-1} \displaystyle \sum_{i=0}^{n-2} \eta_{j} (i) \eta_{j}
(i+1)\right]^{2} \right\rbrace}{\frac{1}{n}} = K_{j,2}. \label{A3:28}
\end{eqnarray}

Equations (\ref{A3:27})-(\ref{A3:28}) imply, for $n$
sufficiently large,
\begin{eqnarray}
&& {\rm Var} \left\lbrace \frac{1}{n} \displaystyle \sum_{i=0}^{n-1} \eta_{j}^{2} (i)\right\rbrace \leq \frac{\widetilde{K}_{j,1}}{n}, \label{A3:27a} \\
&& {\rm Var} \left\lbrace \frac{1}{n-1} \displaystyle \sum_{i=0}^{n-2}
\eta_{j} (i) \eta_{j} (i+1) \right\rbrace \leq
\frac{\widetilde{K}_{j,2}}{n}, \label{A3:28a}
\end{eqnarray}

\noindent for  certain  positive constants  $\widetilde{K}_{j,1}$
and $\widetilde{K}_{j,2},$    for each $j\geq 1.$   Equivalently, for $n$ sufficiently large,
\begin{eqnarray}
{\rm E} \left\lbrace \left(1 - \frac{1}{n} \displaystyle \sum_{i=0}^{n-1} \eta_{j}^{2} (i) \right)^2 \right\rbrace &\leq& \frac{\widetilde{K}_{j,1}}{n}, \label{A3:47b} \\
{\rm E} \left\lbrace \left(\rho_j - \frac{1}{n-1} \displaystyle \sum_{i=0}^{n-1}
\eta_{j} (i)\eta_{j} (i+1)  \right)^2 \right\rbrace &\leq&
\frac{\widetilde{K}_{j,2}}{n}, \label{A3:48b}
\end{eqnarray}

 The following
assumption is now considered.

\bigskip

\noindent \textcolor{Aquamarine}{\textbf{Assumption A4.}} We assume that $$S = \displaystyle \sup_{j
\geq 1} \left( \widetilde{K}_{j,1}+ \widetilde{K}_{j,2} \right) <
\infty.$$

\bigskip

\begin{remark}
\label{A3:remark3}

From equation (\ref{A3:25}), applying the Cauchy--Schwarz's inequality, we
obtain, for each $j \geq 1$,
\begin{eqnarray}
\left| \widehat{\rho}_{n,j} \right| &=& \frac{n}{n-1} \left|
\frac{\displaystyle \sum_{i=0}^{n-2} X_{i,j}X_{i+1,j}}{\displaystyle
\sum_{i=0}^{n-1} X_{i,j}^{2}} \right|\leq
\frac{n}{n-1}\frac{\sqrt{\displaystyle \sum_{i=0}^{n-2}
X_{i,j}^{2}\displaystyle
\sum_{i=0}^{n-2}X_{i+1,j}^{2}}}{\displaystyle \sum_{i=0}^{n-1}
X_{i,j}^{2}}
\nonumber\\
&\leq & \frac{n}{n-1}\sqrt{\frac{\displaystyle
\sum_{i=0}^{n-2}X_{i+1,j}^{2}}{\displaystyle \sum_{i=0}^{n-1}
X_{i,j}^{2}}}\leq  \frac{n}{n-1}~a.s. \label{A3:inef}
\end{eqnarray}

\end{remark}

\newpage

\textcolor{Crimson}{\subsection{Convergence in $\mathcal{L}_{\mathcal{S}(H)}^{2} \left(\Omega, \mathcal{A}, \mathcal{P} \right)$} \label{A3:sec:31}}

Next, the  convergence of $\widehat{\rho}_{k_n}$
to $\rho,$ in the space $\mathcal{L}_{\mathcal{S} \left( H
\right)}^{2} \left( \Omega, \mathcal{A}, \mathcal{P} \right),$ is
derived under the setting of conditions formulated in the previous
sections.

\bigskip
\begin{proposition}
\label{A3:proposition2}

 Let $X= \left\lbrace X_n, \ n \in \mathbb{Z} \right\rbrace$
be a zero--mean standard ARH(1) process. Under \textcolor{Aquamarine}{\textbf{Assumptions A1--A4}},  the following limit holds:
\begin{equation}
\displaystyle \lim_{n \to \infty} \left\| \rho -
\widehat{\rho}_{k_n} \right\|_{\mathcal{L}_{\mathcal{S} \left(H
\right)}^{2} \left( \Omega, \mathcal{A}, \mathcal{P} \right)}^{2} =
0. \label{A3:34}
\end{equation}

Specifically,
\begin{equation}
\left\| \rho - \widehat{\rho}_{k_n}
\right\|_{\mathcal{L}_{\mathcal{S} \left(H \right)}^{2} \left(
\Omega, \mathcal{A}, \mathcal{P} \right)}^{2} \leq g(n), \quad
\mbox{with}\quad g(n) = \mathcal{O} \left( \frac{1}{C_{k_n}^{2} n}
\right),\quad n\rightarrow \infty. \label{A3:33}
\end{equation}
\end{proposition}

\bigskip

\begin{remark}
\label{A3:remark:nuevo_}
\textit{ \cite[Corollary 4.3]{Bosq00} can be applied to obtain weak convergence results,  in terms of weak expectation, using the empirical eigenvectors  . See definition of weak expectation at  the beginning of \cite[Section 1.3, p. 27]{Bosq00}).}
\end{remark}

\bigskip

\begin{proof}
For each $j \geq 1,$ the following almost surely inequality is
satisfied:
\begin{eqnarray}
 \left| \rho_j - \widehat{\rho}_{n,j} \right| &=& \left|
\frac{D_j}{C_j} - \frac{\widehat{D}_{n,j}}{\widehat{C}_{n,j}}
\right| = \left| \frac{D_j \widehat{C}_{n,j} - \widehat{D}_{n,j}C_j}{C_j \widehat{C}_{n,j}} \right| \nonumber \\
&=& \left| \frac{D_j \widehat{C}_{n,j} - \widehat{D}_{n,j}C_j + \widehat{C}_{n,j}\widehat{D}_{n,j} - \widehat{C}_{n,j}\widehat{D}_{n,j}}{C_j \widehat{C}_{n,j}} \right|  \nonumber \\
&=&  \left| \frac{D_j  - \widehat{D}_{n,j}}{C_j} +
\frac{\widehat{C}_{n,j} - C_j}{C_j}
\frac{\widehat{D}_{n,j}}{\widehat{C}_{n,j}}  \right| \leq \frac{1}{C_j} \left(\left| \widehat{\rho}_{n,j} \ \right|
\left| C_j - \widehat{C}_{n,j} \right| + \left| D_j  -
\widehat{D}_{n,j} \right|\right). \nonumber 
\end{eqnarray}

Thus, under \textcolor{Aquamarine}{\textbf{Assumptions A1--A2}}, from equation (\ref{A3:inef}),
for each $j\geq 1,$
\begin{eqnarray}
\left( \rho_j - \widehat{\rho}_{n,j} \right)^{2} &\leq&
\frac{1}{C_{j}^{2}} \left(\left| \widehat{\rho}_{n,j} \ \right|
\left| C_j  - \widehat{C}_{n,j} \right| + \left| D_j  -
\widehat{D}_{n,j} \right|\right)^{2} \nonumber\\ 
&\leq&
\frac{2}{C_{j}^{2}} \left(\left( \widehat{\rho}_{n,j} \right)^2
\left( C_j  - \widehat{C}_{n,j} \right)^2 +
\left( D_j  - \widehat{D}_{n,j} \right)^2 \right) \nonumber \\
&\leq&  \frac{2}{C_{j}^{2}} \left(\left(\frac{n}{n-1} \right)^2
\left( C_j  - \widehat{C}_{n,j} \right)^2 + \left( D_j  -
\widehat{D}_{n,j} \right)^2 \right)~a.s., \nonumber 
\end{eqnarray}

\noindent which implies, for each $j \geq 1$,
\begin{equation}
{\rm E} \left\lbrace \left( \rho_j - \widehat{\rho}_{n,j} \right)^2 \right\rbrace \leq
\frac{2}{C_{j}^{2}} \left( \left(\frac{n}{n-1} \right)^2
{\rm E} \left\lbrace \left( C_j - \widehat{C}_{n,j} \right)^2\right\rbrace  +
{\rm E} \left\lbrace \left( D_j - \widehat{D}_{n,j} \right)^2\right \rbrace \right).
\label{A3:39}
\end{equation}

 Under \textcolor{Aquamarine}{\textbf{Assumption A2}}, from
equations (\ref{A3:24}) and (\ref{A3:39}),

\begin{eqnarray}
\Vert \rho - \widehat{\rho}_{k_n} \Vert_{\mathcal{L}_{\mathcal{S} \left(H \right)}^{2} \left( \Omega, \mathcal{A}, \mathcal{P} \right)}^{2} &=&
{\rm E} \left\lbrace \left\| \rho - \widehat{\rho}_{k_n} \right\|_{\mathcal{S} \left(H \right)}^{2} \right\rbrace = \displaystyle \sum_{j=1}^{k_n} {\rm E} \left\lbrace \left( \rho_j - \widehat{\rho}_{n,j} \right)^2 \right\rbrace + \displaystyle \sum_{j=k_n + 1}^{\infty} {\rm E} \left\lbrace \rho_{j}^{2} \right\rbrace \nonumber \\
&\leq& \displaystyle \sum_{j=1}^{k_n}\frac{2}{C_{j}^{2}} \left(
\left(\frac{n}{n-1} \right)^2 {\rm E} \left\lbrace \left(C_j - \widehat{C}_{n,j}
\right)^2\right\rbrace\right.
\nonumber \\
 & & \left.+{\rm E} \left\lbrace \left( D_j - \widehat{D}_{n,j}\right)^2 \right\rbrace \right) + \displaystyle \sum_{j=k_n + 1}^{\infty}\rho_{j}^{2} \nonumber \\
&\leq&  \frac{2}{C_{k_n}^{2}} \displaystyle \sum_{j=1}^{k_n}  \left(\frac{n}{n-1} \right)^2 \left({\rm E} \left\lbrace \left(C_j - \widehat{C}_{n,j} \right)^2\right\rbrace
\right.\nonumber \\
& &\left.+{\rm E} \left\lbrace \left( D_j - \widehat{D}_{n,j}\right)^2 \right\rbrace\right)  + \displaystyle \sum_{j=k_n + 1}^{\infty}\rho_{j}^{2} \nonumber \\
&\leq&  \frac{2 \left(\frac{n}{n-1} \right)^2}{C_{k_n}^{2}}
\displaystyle \sum_{j=1}^{k_n}  \left({\rm E} \left\lbrace \left(C_j -
\widehat{C}_{n,j} \right)^2\right\rbrace +  {\rm E} \left\lbrace \left( D_j - \widehat{D}_{n,j}\right)^2
\right\rbrace\right) \nonumber \\
&+& \displaystyle \sum_{j=k_n + 1}^{\infty}\rho_{j}^{2}.
 \label{A3:43}
\end{eqnarray}

Furthermore, from  (\ref{A3:15}) and (\ref{A3:25}), for each $j \geq 1$,
\begin{eqnarray}
\widehat{C}_{n,j} &=& \frac{1}{n} \displaystyle \sum_{i=0}^{n-1} X_{i,j}^{2} = \frac{1}{n} \displaystyle \sum_{i=0}^{n-1} C_{j} \eta_{j}^{2} (i), \label{A3:40}\\
\widehat{D}_{n,j} &=& \frac{1}{n-1} \displaystyle \sum_{i=0}^{n-2}
X_{i,j}X_{i+1,j} = \frac{1}{n-1} \displaystyle \sum_{i=0}^{n-2}
C_{j} \eta_j (i)\eta_j (i+1), \label{A3:41}
\end{eqnarray}

\noindent where, considering equation  (\ref{A3:eqproyar1}),
\begin{eqnarray}
D_j &=&  {\rm E} \left\lbrace X_{n,j} X_{n+1,j} \right\rbrace = C_j {\rm E} \left\lbrace \eta_j
(n)\eta_j (n+1)\right\rbrace = C_j \rho_j, \label{A3:42}
\end{eqnarray}

\noindent  for each $j \geq 1.$  Equations (\ref{A3:43})--(\ref{A3:42}) then lead to
\begin{eqnarray}
\Vert \rho - \widehat{\rho}_{k_n} \Vert_{\mathcal{L}_{\mathcal{S} \left(H \right)}^{2} \left( \Omega, \mathcal{A}, \mathcal{P} \right)}^{2} &\leq & \frac{2 \left(\frac{n}{n-1} \right)^2}{C_{k_n}^{2}}
 \displaystyle \sum_{j=1}^{k_n} C_{j}^{2} \left({\rm E} \left\lbrace \left( 1 - \frac{1}{n} \displaystyle \sum_{i=0}^{n-1} \eta_{j}^{2} \left(i
 \right)\right)^2\right\rbrace\right.
\nonumber \\
 & & \left. + {\rm E} \left\lbrace \left(\rho_j - \frac{1}{n-1}\displaystyle \sum_{i=0}^{n-2}  \eta_j (i+1) \eta_j (i) \right)^2 \right\rbrace \right) \nonumber \\
&+& \displaystyle \sum_{j=k_n + 1}^{\infty} \rho_{j}^{2}. \nonumber 
\end{eqnarray}

For each $j \geq 1,$ and for $n$ sufficiently large, considering  equations
(\ref{A3:47b})--(\ref{A3:48b}),  under \textcolor{Aquamarine}{\textbf{Assumption A4}},
\begin{eqnarray}
{\rm E} \left\lbrace \left\| \rho - \widehat{\rho}_{k_n} \right\| _{\mathcal{S} \left(H
\right)}^{2} \right\rbrace &\leq & \frac{2 \left(\frac{n}{n-1}
\right)^2}{C_{k_n}^{2}} \displaystyle
 \sum_{j=1}^{k_n} C_{j}^{2} \left(\frac{\widetilde{K}_{j,1} + \widetilde{K}_{j,2}}{n}\right) + \displaystyle \sum_{j=k_n + 1}^{\infty} \rho_{j}^{2}
\nonumber \\
 &\leq & \frac{2 S \left(\frac{n}{n-1} \right)^2}{
 C_{k_n}^2 n} \displaystyle \sum_{j=1}^{k_n} C_{j}^{2} + \displaystyle \sum_{j=k_n + 1}^{\infty} \rho_{j}^{2}. \nonumber \\
 \label{A3:49}
\end{eqnarray}

From  the trace property of operator $C,$ 

\begin{equation}
\displaystyle \lim_{n\rightarrow \infty}
\displaystyle\sum_{j=1}^{k_n} C_{j}^{2}=
\displaystyle\sum_{j=1}^{\infty}C_{j}^{2}<\infty,\label{A3:trC}
\end{equation}  

\noindent and from the
Hilbert--Schmidt property of $\rho,$   

\begin{equation} \displaystyle
\lim_{n\rightarrow \infty} \displaystyle \sum_{j=k_n + 1}^{\infty}
\rho_{j}^{2}=0.\label{A3:HSrho}
\end{equation} 

 Thus, in view of equations (\ref{A3:49})--(\ref{A3:HSrho}),
\begin{eqnarray}
\left\|  \rho - \widehat{\rho}_{k_n} \right\|
_{\mathcal{L}_{\mathcal{S} \left(H \right)}^{2} \left( \Omega,
\mathcal{A}, \mathcal{P} \right)}^{2} &=& {\rm E} \left\lbrace \left\|  \rho -
\widehat{\rho}_{k_n} \right\| _{\mathcal{S} \left(H \right)}^{2} \right\rbrace 
\leq g(n) = \mathcal{O} \left(\frac{1}{C_{k_n}^2 n} \right), \
n\rightarrow \infty, \nonumber\\ \label{A3:50}
\end{eqnarray}

\noindent where \begin{equation}g(n) = \frac{2 S \left(\frac{n}{n-1}
\right)^2}{C_{k_n}^2 n} \displaystyle \sum_{j=1}^{k_n} C_{j}^{2} +
\displaystyle \sum_{j=k_n + 1}^{\infty} \rho_{j}^{2}.\label{A3:eqdefg}
\end{equation}

Under \textcolor{Aquamarine}{\textbf{Assumption A3}}, equations (\ref{A3:50})--(\ref{A3:eqdefg}) imply
$$\lim_{n\rightarrow \infty}\Vert \rho - \widehat{\rho}_{k_n} \Vert_{\mathcal{L}_{\mathcal{S} \left(H \right)}^{2} \left( \Omega, \mathcal{A}, \mathcal{P} \right)}^{2}=0,$$
as we wanted to prove. 

\hfill \hfill \textcolor{Aquamarine}{$\blacksquare$}

\end{proof}

\bigskip

Note that consistency of $\widehat{\rho}_{k_n}$ in the space
$\mathcal{S} \left(H \right)$ directly follows from equation
(\ref{A3:34}) in \textcolor{Crimson}{Proposition} \ref{A3:proposition2}.

\bigskip

\begin{corollary}
\label{A3:cor1} Let $X= \left\lbrace X_n,\ n \in \mathbb{Z} \right \rbrace$ be a
zero--mean standard ARH(1) process. Under \textcolor{Aquamarine}{\textbf{Assumptions A1--A4}},
as long as $n\rightarrow \infty,$
\begin{equation}
\left\|  \rho - \widehat{\rho}_{k_n} \right\| _{\mathcal{S} \left(H
\right)} \to^{p} 0, \nonumber 
\end{equation}
\noindent where, as usual, $\to^{p}$ denotes the convergence in probability.
\end{corollary}

%
%

\textcolor{Crimson}{\subsection{Consistency of the ARH(1)  plug--in  predictor.} \label{A3:sec:32}}

Let us consider $\mathcal{L} \left( H \right)$ the space  of
bounded linear operators on $H,$  with the norm
\begin{equation}
\left\| \mathcal{A} \right\|_{\mathcal{L} \left( H \right)} =
\displaystyle \sup_{x \in H} \frac{\left\| \mathcal{A} \left( x
\right) \right\|_H}{\left\| x \right\|_H}, \nonumber 
\end{equation}
\noindent for every  $\mathcal{A}\in \mathcal{L} \left( H \right).$
In particular,  for each  $x \in H,$
\begin{equation}
\left\| \mathcal{A} \left( x\right) \right\|_H \leq \left\|
\mathcal{A} \right\|_{\mathcal{L} \left( H \right)}\left\| x
\right\|_H. \label{A3:59}
\end{equation}

In the following, we denote by  
\begin{equation}\widehat{X}_{n} =
\widehat{\rho}_{k_n} \left( X_{n-1} \right)\label{A3:ARHpred}
\end{equation}

\noindent  as usual, the ARH(1) plug--in
predictor of $X_n,$ as an estimator of the conditional expectation \linebreak
${\rm E} \left\lbrace X_{n}|X_{n-1} \right\rbrace= \rho \left( X_{n-1} \right)$.  The following
proposition provides the consistency of $\widehat{X}_{n} =
\widehat{\rho}_{k_n} \left( X_{n-1} \right)$ in $H$.

 \bigskip

\begin{proposition}
\label{A3:cor2} 
\textit{Let $X = \left\lbrace X_n, \ n \in \mathbb{Z} \right\rbrace$ be a
zero--mean standard ARH(1) process. Under \textcolor{Aquamarine}{\textbf{Assumptions A1--A4}},
\begin{equation}
\lim_{n\rightarrow \infty }{\rm E} \left\lbrace \left\| \left(\rho - \widehat{\rho}_{k_n} \right) \left( X_{n-1} \right) \right\|_H \right\rbrace
=0. \nonumber 
\end{equation}
Specifically,
\begin{equation}
{\rm E} \left\lbrace \left\| \left(\rho - \widehat{\rho}_{k_n} \right) \left( X_{n-1} \right)   \right\|_H \right\rbrace \leq h \left(n
\right),\quad h \left( n \right) = \mathcal{O}
\left(\frac{1}{C_{k_n} \sqrt{n}} \right),\quad n\rightarrow \infty. \nonumber 
\end{equation}
In particular,
\begin{equation}
\left\| \left(\rho - \widehat{\rho}_{k_n} \right) \left( X_{n-1} \right)  \right\|_H \to^{p} 0, \nonumber 
\end{equation}
where, as usual, $\to^{p}$ denotes the convergence in probability.}
\end{proposition}

\begin{proof}

From (\ref{A3:59}) and \textcolor{Crimson}{Proposition} \ref{A3:proposition2}, for $n$ sufficiently large, the following almost surely inequality holds:
\begin{equation}
\left\| \rho \left(X_{n-1} \right) - \widehat{X}_{n}  \right\|_H\leq \left\| \rho -
\widehat{\rho}_{k_n} \right\|_{\mathcal{L} \left( H
\right)}\left\|X_{n-1}\right\|_{H}, \nonumber 
\end{equation}

\noindent where, as given in equation (\ref{A3:ARHpred}), $\widehat{X}_{n}=\widehat{\rho}_{k_n} \left( X_{n-1} \right).$
Thus,
\begin{equation}
{\rm E} \left\lbrace  \left\| \rho \left(X_{n-1} \right) - \widehat{X}_{n}   \right\|_H \right\rbrace \leq {\rm E} \left\lbrace  \left\|
\rho - \widehat{\rho}_{k_n} \right\|_{\mathcal{L} \left( H \right)}
\left\|X_{n-1}\right\|_{H} \right\rbrace. \label{A3:60}
\end{equation}

From the Cauchy-Schwarz's inequality, keeping in mind that, for a
Hilbert--Schmidt operator $\mathcal{K},$ it always holds that
 $\|\mathcal{K}\|_{\mathcal{L} \left( H \right)}\leq
\|\mathcal{K}\|_{\mathcal{S} \left( H \right)},$ we have from
equation (\ref{A3:60}),
\begin{eqnarray}
{\rm E} \left\lbrace \left\| X_n - \widehat{X}_{n}  \right\|_H \right\rbrace & \leq & \sqrt{ {\rm E} \left\lbrace \left\|
\rho - \widehat{\rho}_{k_n} \right\|_{\mathcal{L} \left( H
\right)}^{2} \right\rbrace} \sqrt{{\rm E} \left\lbrace \left\|X_{n-1}\right\|_{H}^{2}\right\rbrace } \nonumber\\ 
& \leq & \sqrt{{\rm E} \left\lbrace\left\| \rho - \widehat{\rho}_{k_n}
\right\|_{\mathcal{S}
 \left( H \right)}^{2} \right\rbrace} \sqrt{{\rm E} \left\lbrace \left\|X_{n-1}\right\|_{H}^{2} \right\rbrace }\nonumber \\
 &=& \sqrt{{\rm E} \left\lbrace \left\| \rho - \widehat{\rho}_{k_n} \right\|_{\mathcal{S}
 \left( H \right)}^{2}\right\rbrace}\sigma_{X}, \label{A3:64}
\end{eqnarray}

\noindent where, as before,  $\sigma_{X}^{2} = {\rm E} \left\lbrace
\left\|X_{n-1}\right\|_{H}^{2} \right\rbrace  = \displaystyle \sum_{j=1}^{\infty}
C_j < \infty, \quad n \in \mathbb{Z}$ (see equation (\ref{A3:20})).

Since from \textcolor{Crimson}{Proposition} \ref{A3:proposition2} (see equation
(\ref{A3:33})),$$\left\| \rho - \widehat{\rho}_{k_n}
\right\|_{\mathcal{L}_{\mathcal{S} \left(H \right)}^{2} \left(
\Omega, \mathcal{A}, \mathcal{P} \right)}^{2} \leq g(n), \quad
\mbox{with}\quad g(n) = \mathcal{O} \left( \frac{1}{C_{k_n}^{2} n}
\right),\quad n\rightarrow \infty,$$  from equation
(\ref{A3:64}), we obtain,
\begin{equation}
{\rm E} \left\lbrace \left\|  \left( \rho - \widehat{\rho}_{k_n} \right) \left(X_{n-1} \right)  \right\|_H \right\rbrace  \leq h \left( n
\right), \nonumber
\end{equation}

\noindent where $h \left( n \right) = \sigma_{X}  \sqrt{g \left( n
\right)}$, with $g \left( n \right)$ being given in (\ref{A3:eqdefg}).
In particular, under \textcolor{Aquamarine}{\textbf{Assumption A3}},
\begin{equation}
\displaystyle \lim_{n \to \infty} {\rm E} \left\lbrace \left\|  \left( \rho - \widehat{\rho}_{k_n} \right) \left(X_{n-1} \right) 
\right\|_H \right\rbrace = 0, \nonumber 
\end{equation}

\noindent which implies that
\begin{equation}
\left\|  \left( \rho - \widehat{\rho}_{k_n} \right) \left(X_{n-1} \right)    \right\|_H = \left\| \rho  \left(X_{n-1} \right) - \widehat{X}_{n}   \right\|_H \to^{p} 0,\quad n\rightarrow
\infty. \nonumber 
\end{equation}

\hfill \hfill \textcolor{Aquamarine}{$\blacksquare$}

\end{proof}

%
%

\textcolor{Crimson}{\section{The Gaussian case}
\label{A3:gc}} 

In this section, we prove that, in the Gaussian ARH(1) context, \textcolor{Aquamarine}{\textbf{Assumptions A1--A2}} and \textcolor{Aquamarine}{\textbf{A4}} also hold. From equation (\ref{A3:22}), for $n \geq 1$,

\begin{equation}
{\rm E} \left\lbrace \frac{\displaystyle \sum_{i=0}^{n-1}\eta_{j}^{2}(i)}{n}\right\rbrace=1. \nonumber 
\end{equation}

Furthermore, for each $j \geq 1$ and $n \geq 2$, the $n\times 1$ random vector $\boldsymbol{\eta
}_{j}^{T}=\left(\eta_{j}(0),\dots,\eta_{j}(n-1)\right)$ 
follows a Multivariate Normal
distribution with null mean vector, and covariance matrix
\begin{equation}
\boldsymbol{\Sigma} = \left(\begin{array}{cccccc}
1 & \rho_{j} & 0 & \dots &\dots & 0\\
\rho_{j} & 1 & \rho_{j}& 0 &\dots & 0\\
0 & \rho_{j} & 1 & \rho_{j} & \dots &0\\
\vdots &\vdots& \vdots & \vdots & \vdots &\vdots \\
0 & 0 & 0 & 0 & \rho_{j} & 1\\
\end{array}\right)_{n\times n}.\label{A3:eqmatex1}
\end{equation}

 It is well--known (see, for
example, \cite{Gurland56})  that the  variance of a quadratic form
defined from a multivariate Gaussian vector $\mathbf{y}\sim
N(\boldsymbol{\mu},\boldsymbol{\Lambda} ),$ and a symmetric matrix $\mathbf{Q}$
is given by:
\begin{eqnarray}
{\rm Var} \left\lbrace \mathbf{y}^{T}\boldsymbol{Q} \mathbf{y}\right\rbrace &=& 2
{\rm Tr} \left(\boldsymbol{Q} \boldsymbol{\Lambda
Q\Lambda}\right)+4\boldsymbol{\mu}^{T}\boldsymbol{Q\Lambda
Q}\boldsymbol{\mu}.\label{A3:mvqf2}
\end{eqnarray}

For each $j\geq 1,$ applying equation (\ref{A3:mvqf2}), with
$\mathbf{y}=\boldsymbol{\eta }_{j},$ $\boldsymbol{\Lambda =\Sigma}$ in (\ref{A3:eqmatex1}),
and $\boldsymbol{Q=Id}_{n},$ the $n\times n$ identity matrix, keeping in
mind ${\rm E} \left\lbrace \eta_j (i) \eta_j (i+1) \right\rbrace = \rho_j$, for every $i \in \mathbb{Z}$,

\begin{eqnarray}
{\rm Var} \left\lbrace \boldsymbol{\eta }_{j}^{T} \boldsymbol{Id}_{n}\boldsymbol{\eta }_{j}\right\rbrace &=&
{\rm Var} \left\lbrace \sum_{i=0}^{n-1}\eta_{j}^{2}(i)\right\rbrace
= 2
{\rm Tr} \left(\boldsymbol{\Sigma \Sigma}
\right)=2\left(n+2(n-1)\rho_{j}^{2}\right). \nonumber \\\label{A3:mvqf4}
\end{eqnarray}

Furthermore, from equation (\ref{A3:mvqf4}), for each $j\geq 1,$
\begin{eqnarray}
{\rm Var} \left\lbrace \frac{\displaystyle \sum_{i=0}^{n-1}\eta_{j}^{2}(i)}{n}\right\rbrace=
\frac{2}{n^{2}}\left(n+2(n-1)\rho_{j}^{2}\right)=\frac{2}{n}+4\left(\frac{1}{n}-\frac{1}{n^{2}}\right)\rho_{j}^{2}.
\label{A3:expectation2}
\end{eqnarray}

We then obtain, from equation (\ref{A3:expectation2}), 

\begin{eqnarray}
\lim_{n\rightarrow
\infty}{\rm Var} \left\lbrace \frac{\displaystyle \sum_{i=0}^{n-1}\eta_{j}^{2}(i)}{n}\right\rbrace &=&
\lim_{n\rightarrow
\infty}{\rm E} \left\lbrace \left(1-\frac{\displaystyle \sum_{i=0}^{n-1}\eta_{j}^{2}(i)}{n}\right)^{2}\right\rbrace
\nonumber\\
&=&\lim_{n\rightarrow
\infty}\frac{2}{n}+4\left(\frac{1}{n}-\frac{1}{n^{2}}\right)\rho_{j}^{2}=0.\label{A3:zerolimitbb}
\end{eqnarray}

Equation (\ref{A3:zerolimitbb}) leads to

\begin{equation}
\lim_{n\rightarrow \infty}\frac{{\rm Var}\left\lbrace \frac{\displaystyle \sum_{i=0}^{n-1}\eta_{j}^{2}(i)}{n}\right\rbrace}{\frac{1}{n}}=2+4\rho_{j}^{2}. \nonumber 
\end{equation}

Hence, for each $j\geq 1,$ $K_{j,1}$ in equation (\ref{A3:27})  is
given by 
\begin{equation}
K_{j,1}= 2+4\rho_{j}^{2},  \nonumber 
\end{equation}

\noindent and, from equation (\ref{A3:expectation2}), 
$${\rm Var} \left\lbrace \frac{\displaystyle\sum_{i=0}^{n-1}\eta_{j}^{2}(i)}{n}\right\rbrace \leq 2 +  4 \left(\frac{1}{n} - \frac{1}{n^2} \right)
\rho_{j}^{2} \leq 2 +  4 \rho_{j}^{2} \leq 6.$$ 

Thus, for every $j\geq 1,$ $\widetilde{K}_{j,1}$ in equation
(\ref{A3:27a}) satisfies
\begin{equation}
\widetilde{K}_{j,1}\leq 6. \nonumber 
\end{equation}

\bigskip
\begin{remark}
\textit{Note that, from \textcolor{Crimson}{Lemma} \ref{A3:proposition1},  for each $j\geq 1$ and $i \in \mathbb{Z}$,
\begin{equation} 
{\rm E} \left\lbrace \widetilde{\eta}_{j}^{4}(i)\right\rbrace=3. \nonumber 
\end{equation}
 Thus,  the assumption considered in \textcolor{Crimson}{Remark} \ref{A3:remark1} holds, and    for each $j\geq 1,$ the AR(1) process \linebreak $\eta_{j}=\left\lbrace \eta_{j}(n),\ n\in \mathbb{Z} \right\rbrace$ is
  ergodic for all second--order moments, in the mean--square sense; see \cite[pp. 192--193]{Hamilton94}.}
\end{remark}

\bigskip

For  $n\geq 2,$ and for each $j\geq 1,$ we are now going to compute
$K_{j,2}$ in (\ref{A3:28}). The $(n-1)\times 1$ random vectors  $$\boldsymbol{\eta
}_{j}^{\star}=\left(\eta_{j}(0),\dots,\eta_{j}(n-2)\right)^{T}, \quad \boldsymbol{\eta }_{j}^{\star\star
}=\left(\eta_{j}(1),\dots,\eta_{j}(n-1)\right)^{T}$$ are
multivariate Normal distributed, with null   mean vector,   and  covariance matrix
\begin{equation}
\boldsymbol{\widetilde{\Sigma}}=
\left(\begin{array}{cccccc}
1 & \rho_{j} & 0 & \dots &\dots & 0\\
\rho_{j} & 1 & \rho_{j}& 0 &\dots & 0\\
0 & \rho_{j} & 1 & \rho_{j} & \dots &0\\
\vdots &\vdots& \vdots & \vdots & \vdots &\vdots \\
0 &\dots& 0 & 0 & \rho_{j} & 1\\
\end{array}\right)_{(n-1)\times (n-1)}.
\label{A3:eqcov2}
\end{equation}

From equation (\ref{A3:24a}), for each $j\geq 1,$

\begin{equation}
{ \rm E} \left\lbrace\sum_{i=0}^{n-2}\eta_{j}(i)\eta_{j}(i+1)\right\rbrace=\sum_{i=0}^{n-2}\rho_{j}=(n-1)\rho_{j}= {\rm Tr} \left( { \rm E} \left\lbrace \boldsymbol{\eta }_{j}^{\star}[\boldsymbol{\eta
}_{j}^{\star\star }]^{T}\right\rbrace\right),\label{A3:eqcovmbb}
\end{equation}

\noindent where
\begin{equation}
{ \rm E} \left\lbrace \boldsymbol{\eta }_{j}^{\star }[\boldsymbol{\eta
}_{j}^{\star \star}]^{T}\right\rbrace = { \rm E} \left\lbrace \boldsymbol{\eta }_{j}^{\star
}\otimes \boldsymbol{\eta }_{j}^{\star
\star}\right\rbrace =\rho_{j} \boldsymbol{Id}_{n-1},\label{A3:eqcovm}
\end{equation}

\noindent with, as before, $\boldsymbol{Id}_{n-1}$ denoting the
$(n-1)\times (n-1)$ identity matrix.

However, the
variance of $$\displaystyle \sum_{i=0}^{n-2}\eta_{j}(i)\eta_{j}(i+1)$$ depends
greatly on the distribution of $\boldsymbol{\eta }_{j}^{\star }$ and
$\boldsymbol{\eta }_{j}^{\star \star}.$  In the Gaussian case,
keeping in mind that  $$\boldsymbol{\eta
}_{j}^{\star}=\left(\eta_{j}(0),\dots,\eta_{j}(n-2)\right)^{T}, \quad \boldsymbol{\eta }_{j}^{\star\star
}=\left(\eta_{j}(1),\dots,\eta_{j}(n-1)\right)^{T}$$  are
zero--mean multivariate Normal distributed vectors with covariance matrix
$\boldsymbol{\widetilde{\Sigma}}$ given in (\ref{A3:eqcov2}), and having
cross--covariance matrix in (\ref{A3:eqcovm}), we can  compute the
variance of $\displaystyle \sum_{i=0}^{n-2}\eta_{j}(i)\eta_{j}(i+1),$ from  (\ref{A3:eqcovmbb})--(\ref{A3:eqcovm}), as follows. First,

\begin{eqnarray}
{ \rm Var} \left\lbrace [\boldsymbol{\eta
}_{j}^{\star}]^{T} \boldsymbol{Id}_{n-1}\boldsymbol{\eta
}_{j}^{\star\star }\right\rbrace &=& { \rm E} \left\lbrace [\boldsymbol{\eta
}_{j}^{\star}]^{T} \boldsymbol{Id}_{n-1}\boldsymbol{\eta
}_{j}^{\star\star }[\boldsymbol{\eta
}_{j}^{\star}]^{T} \boldsymbol{Id}_{n-1}\boldsymbol{\eta
}_{j}^{\star\star }\right\rbrace \nonumber\\
& -&\left( { \rm E} \left\lbrace [\boldsymbol{\eta }_{j}^{\star}]^{T}\boldsymbol{Id}_{n-1}\boldsymbol{\eta }_{j}^{\star\star
}\right\rbrace]\right)^{2}. \nonumber
\end{eqnarray}

This can be rewritten as
\begin{equation}
\sum_{i=0}^{n-2}\sum_{p=0}^{n-2}
{ \rm E} \left\lbrace \eta_{j}(i)\eta_{j}(i+1)\eta_{j}(p)\eta_{j}(p+1) \right\rbrace
-\left({ \rm E} \left\lbrace [\boldsymbol{\eta }_{j}^{\star}]^{T} \boldsymbol{Id}_{n-1}\boldsymbol{\eta }_{j}^{\star\star
}\right\rbrace\right)^{2}, \nonumber
\end{equation}
\noindent which is equal to
\begin{eqnarray}
\sum_{i=0}^{n-2} { \rm E} \left\lbrace \eta_{j}(i)\eta_{j}(i+1)\right\rbrace \sum_{p=0}^{n-2} { \rm E} \left\lbrace \eta_{j}(p)\eta_{j}(p+1)\right\rbrace
&+&
\sum_{i=0}^{n-2}\sum_{p=0}^{n-2}{ \rm E} \left\lbrace \eta_{j}(i)\eta_{j}(p)\right\rbrace { \rm E} \left\lbrace \eta_{j}(i+1)\eta_{j}(p+1)\right\rbrace
\nonumber\\
&+&\sum_{i=0}^{n-2}\sum_{p=0}^{n-2} { \rm E} \left\lbrace \eta_{j}(i)\eta_{j}(p+1)\right\rbrace { \rm E} \left\lbrace \eta_{j}(i+1)\eta_{j}(p)\right\rbrace \nonumber \\
&-& \left( { \rm E} \left\lbrace[\boldsymbol{\eta }_{j}^{\star}]^{T}\boldsymbol{Id}_{n-1}\boldsymbol{\eta }_{j}^{\star\star }\right\rbrace\right)^{2}. \nonumber
\end{eqnarray}

This then reduces to
\begin{eqnarray}
\left[{\rm Tr} \left({ \rm E} \left\lbrace \boldsymbol{\eta
}_{j}^{\star}\otimes \boldsymbol{\eta }_{j}^{\star\star
}\right\rbrace\right)\right]^{2}&+& {\rm Tr} \left(\boldsymbol{\widetilde{\Sigma}} \boldsymbol{\widetilde{\Sigma}}\right)
\nonumber\\
&+&
{\rm Tr} \left({ \rm E} \left\lbrace \boldsymbol{\eta }_{j}^{\star}\otimes
\boldsymbol{\eta }_{j}^{\star\star } \right\rbrace \left[{ \rm E} \left\lbrace \boldsymbol{\eta
}_{j}^{\star}\otimes \boldsymbol{\eta }_{j}^{\star\star
}\right\rbrace \right]^{T}\right)-\left[{\rm Tr}\left({ \rm E} \left\lbrace \boldsymbol{\eta }_{j}^{\star}\otimes
\boldsymbol{\eta }_{j}^{\star\star
}\right\rbrace \right)\right]^{2}, \nonumber \\ \label{A3:extversionvv3}
\end{eqnarray}
\noindent which is the same as
\begin{eqnarray}
{\rm Tr}\left(\boldsymbol{\widetilde{\Sigma}} \boldsymbol{\widetilde{\Sigma}}\right)&+&
{\rm Tr}\left({ \rm E} \left\lbrace \boldsymbol{\eta }_{j}^{\star}\otimes
\boldsymbol{\eta }_{j}^{\star\star }\right\rbrace \left[{ \rm E} \left\lbrace \boldsymbol{\eta
}_{j}^{\star}\otimes \boldsymbol{\eta }_{j}^{\star\star
}\right\rbrace \right]^{T}\right) \nonumber \\
&=&(n-1)+2(n-2)\rho_{j}^{2}+(n-1)\rho_{j}^{2}, \nonumber
\end{eqnarray}

\noindent where, from (\ref{A3:eqcovm}), 
$${ \rm E} \left\lbrace \boldsymbol{\eta
}_{j}^{\star}\otimes \boldsymbol{\eta }_{j}^{\star
\star}\right\rbrace \left[{ \rm E} \left\lbrace \boldsymbol{\eta }_{j}^{\star}\otimes
\boldsymbol{\eta }_{j}^{\star
\star}\right\rbrace \right]^{T}=\left(\begin{array}{ccccc} \rho_{j}^{2}& 0 &\dots & \dots & 0\\
0 &\rho_{j}^{2}&0 &\dots & 0\\
\vdots & \ddots &\ddots  & \vdots & \vdots\\
0 &\dots &\ddots & \ddots& \rho_{j}^{2}\\
\end{array}\right)=\rho_{j}^{2} \boldsymbol{Id}_{n-1}.$$ 

 From
(\ref{A3:extversionvv3}),

\begin{equation}
{ \rm Var} \left\lbrace \frac{\displaystyle \sum_{i=0}^{n-2}\eta_{j}(i)\eta_{j}(i+1)}{n-1}\right\rbrace =
\frac{(n-1)+2(n-2)\rho_{j}^{2}+(n-1)\rho_{j}^{2}}{(n-1)^{2}}.
\label{A3:fvar}
\end{equation}

 Therefore, for each $j\geq 1,$
\begin{equation}
\lim_{n\rightarrow \infty}n{\rm Var}\left\lbrace \frac{\displaystyle \sum_{i=0}^{n-2}\eta_{j}(i)\eta_{j}(i+1)}{n-1}\right\rbrace=1+3\rho_{j}^{2}. \nonumber 
\end{equation}

Thus, for each $j\geq 1,$ $K_{j,2}$ in (\ref{A3:28}) is given by
$K_{j,2}=1+3\rho_{j}^{2}.$ From equation (\ref{A3:fvar}), 
$${\rm Var} \left\lbrace \frac{\displaystyle \sum_{i=0}^{n-2}\eta_{j}(i)\eta_{j}(i+1)}{n-1}\right\rbrace \leq   1 + 3 \rho_{j}^{2} \leq
4.$$ Hence, for every $j\geq 1,$ $\widetilde{K}_{j,2}$ in equation
(\ref{A3:28a}) satisfies
\begin{equation}
\widetilde{K}_{j,2}\leq 4. \nonumber 
\end{equation}

Therefore, the constant
$S$ in \textcolor{Aquamarine}{\textbf{Assumption A4}} is such that $S\leq  6+4=10.$

%
%

\textcolor{Crimson}{\section{Simulation study}
\label{A3:sec:4}}

A simulation study is undertaken to  illustrate the behaviour of the
formulated componentwise estimator of the autocorrelation operator,
and of its associated ARH(1) plug--in predictor for large sample
sizes. The results are reported in \textcolor{Crimson}{Appendix} \ref{A3:secBlss}. In \textcolor{Crimson}{Appendix} \ref{A3:ACS}, a comparative study is developed, from the implementation of the ARH(1) plug--in prediction techniques proposed in \cite{AntoniadisSapatinas03,Besseetal00,Bosq00,Guillas01}. In the subsequent sections, we restrict our attention to the Gaussian case

\textcolor{Crimson}{\subsection{Behaviour of $\widehat{\rho}$ and $\widehat{X}_{n}$ for large sample sizes}
\label{A3:secBlss}}

Let  $(- \Delta)_{(a,b)}$ be the Dirichlet negative Laplacian operator on $(a,b)$  given by

\begin{eqnarray}
(- \Delta)_{(a,b)} \left( f \right) \left( x \right) &=& - \displaystyle  \frac{d^2}{d x^{2}} f \left( x \right),\quad x  \in (a,b)\subset \mathbb{R}, \nonumber \\
f \left(a\right) &=& f(b)=0.  \nonumber 
\end{eqnarray}

The eigenvectors $\left\lbrace \phi_j, \ j \geq 1 \right\rbrace$  and eigenvalues $\left\lbrace \lambda_j \left( (- \Delta)_{(a,b)} \right), \ j \geq 1
\right\rbrace$  of $(- \Delta)_{(a,b)}$ satisfy, for each $j\geq 1$ and for each $x \in (a,b)$,
\begin{eqnarray}
(- \Delta)_{(a,b)} \phi_j \left( x \right) &=& \lambda_j \left( (- \Delta)_{(a,b)} \right) \phi_j \left( x \right) ,\quad  
\phi_j \left( a \right)= \phi_j \left( b \right)=0.
\label{A3:69}
\end{eqnarray}

For each $j \geq 1$ and $x \in \left[a,b \right]$, the solution to equation (\ref{A3:69}) is given by
(see \cite[p. 6]{GrebenkovNguyen13}):
\begin{eqnarray}
\phi_{j} \left( x\right) &=& \sqrt{\frac{2}{b-a}} \displaystyle \sin \left( \frac{\pi j x}{b-a} \right), \quad \forall x \in \left[a, b \right],\quad \lambda_{j} \left( (- \Delta)_{(a,b)} \right) = \frac{\pi^2
j^2}{(b-a)^{2}}. \label{A3:extra__}
\end{eqnarray}

We consider here the operator $C$ defined as
\begin{equation} 
C=\left( (- \Delta)_{(a,b)} \right)^{-2 \left(1 - \gamma_1 \right)}, \quad \gamma_1 \in \left(0,1/2\right). \nonumber 
\end{equation}

From \cite[pp. 119--140]{DautrayLions90}, the eigenvectors   of $C$
coincide with the eigenvectors of $(- \Delta)_{(a,b)},$  and its
eigenvalues $\left\lbrace C_{j}, \ j\geq 1 \right\rbrace$ are given by:
\begin{equation}
C_j = \left[\lambda_j \left( (- \Delta)_{(a,b)} \right)\right] ^{-2 \left(1 - \gamma_1 \right)}= \left[\frac{\pi^2 j^{2}}{(b-a)^{2}}\right]^{-2 \left(1 - \gamma_1 \right)}.\label{A3:eigvC}
\end{equation}

Additionally, considering  
\begin{equation}
\rho=\left[\frac{ (-
\Delta)_{(a,b)} }{{\lambda_1 \left( (- \Delta)_{(a,b)} \right)
 - \epsilon}}\right]^{-\left(1 - \gamma_2 \right)},\quad \gamma_2 \in \left(0,1/2\right),\nonumber 
\end{equation} 
\noindent for certain positive constant $\epsilon <\lambda_1 \left( (- \Delta)_{(a,b)} \right)$ close to zero, $\rho $ is a positive  self--adjoint Hilbert--Schmidt operator, whose eigenvectors coincide with the eigenvectors of
 $\left( - \Delta\right)_{(a,b)},$ and whose eigenvalues $\left\lbrace \rho_{j},~j\geq 1 \right\rbrace$ are such that $\rho_{j}<1,$ for every $j\geq 1,$  and

\begin{equation}
\rho_{j}^{2} = \left[\frac{\lambda_j \left(( - \Delta)_{(a,b)}
\right)}{\lambda_1 \left( (- \Delta)_{(a,b)} \right) -
\epsilon}\right]^{-2\left(1 - \gamma_2 \right)}, \quad
\rho_{j}^{2} \in \left(0, 1 \right),\quad  \gamma_2 \in \left(0,1/2 \right), \label{A3:75}
\end{equation}

\noindent where,  as before,  $\left\lbrace \lambda_j \left( (- \Delta)_{(a,b)}\right), \ j \geq 1 \right\rbrace$ are given in equation (\ref{A3:extra__}).

From
(\ref{A3:25a}), the eigenvalues $\left\lbrace \sigma_{j}^{2}, \ j\geq 1 \right\rbrace$ of
$C_{\varepsilon}$ are then defined, for each $j \geq 1,$ as
\begin{equation}
\sigma_{j}^{2} = C_j \left(1 - \rho_{j}^{2} \right) = [\lambda_j
\left( (- \Delta)_{(a,b)}  \right)]^{-2 \left( 1 - \gamma_1 \right)} -
\frac{\left[\lambda_j \left( (- \Delta)_{(a,b)}  \right) \right]^{-2 \left( 2 - \gamma_1 - \gamma_2 \right)}}{\left[\lambda_1 \left( (-
\Delta)_{(a,b)} \right) - \epsilon \right]^{-2 \left( 1 - \gamma_2 \right)}}. \nonumber 
\end{equation}

 Note that $C_{\varepsilon}$ is in the trace class, since
the trace property of $C,$ and the fact that  $\rho_{j}^{2}< 1,$ for
every $j\geq 1,$ implies $$\displaystyle \sum_{j=1}^{\infty}
\sigma_{j}^{2}=\displaystyle \sum_{j=1}^{\infty}C_j \left(1 -
\rho_{j}^{2} \right) <\displaystyle \sum_{j=1}^{\infty}C_j<
\infty.$$

For this particular example of operator $C,$ we have considered truncation parameter $k_{n}$ of the form
\begin{equation}
k_{n}= n^{1/\alpha}, \label{A3:77}
\end{equation}

\noindent for a suitable $\alpha >0,$  which, in particular, allows verification of  (\ref{A3:26}). From equation (\ref{A3:eigvC}), one has, for $\gamma_1 \in (0, 1/2)$,
\begin{equation}
\sqrt{n}C_{k_n} = \sqrt{n} \left[\lambda_{k_n} \left( - \Delta_{(a,b)}
\right)\right]^{-2\left(1 - \delta_1 \right)} = \sqrt{n}\left( \frac{ \pi k_n
}{b-a} \right)^{-4\left(1 - \delta_1 \right)},\quad ~\delta_1 > 1. \nonumber 
\end{equation}

 From equation (\ref{A3:77}), \textcolor{Aquamarine}{\textbf{Assumption A3}} is then satisfied if
\begin{equation}
1/2- \frac{4 \left(1 - \gamma_1\right)}{\alpha } > 0, \quad \mbox{i.e., \
if}\quad \alpha > 8 \left(1 - \gamma_1 \right) > 4. \label{A3:86}
\end{equation}
\noindent since $\gamma_1 \in  (0,1/2)$. Fix $\gamma_1 = 0.4$ and $\gamma_2 = 9/20$. Then, from equation (\ref{A3:86}),  $\alpha
> 48/10.$ In particular, the values  $\alpha_{1} = 5$ and $\alpha _{2}=6$ have been tested, in Table \ref{A3:tab:Table1} below, for $H=L^{2}((a,b)),$ and  $(a,b)=(0,4),$ where $L^{2}((a,b))$ denotes the space of square integrable functions on $(a,b).$

The
computed empirical truncated functional mean square error
${\rm EMSE}_{\widehat{\rho}_{k_{n}}}$  of the  estimator
$\widehat{\rho}_{k_n}$ of $\rho,$ for a sample size $n$, is given by:
\begin{eqnarray}
{\rm EMSE}_{\widehat{\rho}_{k_{n}}}&=& \frac{1}{N} \displaystyle \sum_{w=1}^{N} \displaystyle \sum_{j=1}^{k_n} \left( \rho_{j} - \widehat{\rho}_{n,j}^{w} \right)^2, \label{A3:92bb}\\
\widehat{\rho}_{n,j}^{w} &=&
\frac{\widehat{D}_{n,j}^{w}}{\widehat{C}_{n,j}^{w}} =
\frac{\frac{1}{n-1} \displaystyle \sum_{i=0}^{n-2}
X_{i,j}^{w}X_{i+1,j}^{w} }{\frac{1}{n} \displaystyle
\sum_{i=0}^{n-1} \left(X_{i,j}^{w} \right)^2}, \label{A3:93bb}
\end{eqnarray}

\noindent where $N$ denotes the number of simulations, and for each
$j=1,\dots, k_n,$ $\widehat{\rho}_{n,j}^{w}$ represents the
estimator of $\rho_{j},$ based on  the $w$--th generation of the
values $X_{0,j}^{w},\dots, X_{n-1,j}^{w},$ with $X_{i,j}^{w}=\left\langle X_{i}^{w},\phi_{j}\right\rangle_{H},$ for \linebreak $w=1,\dots, 700,$ and $i=0,\dots,n-1.$

 For the plug--in predictor
$\widehat{X}_n = \widehat{\rho}_{k_n} \left( X_{n-1} \right),$ we
compute the empirical version
${\rm UB(EMAE)}_{\widehat{X}_n^{k_n}}$ of the derived upper bound
(\ref{A3:64}),  which, for each $n \in \mathbb{Z},$ is given by
\begin{eqnarray}
{\rm UB(EMAE)}_{\widehat{X}_n^{k_n}} &=& \sqrt{\frac{1}{N}
\displaystyle \sum_{w=1}^{N} \displaystyle \sum_{j=1}^{k_n} \left(
\rho_{j} - \widehat{\rho}_{n,j}^{w} \right)^2
\widehat{{ \rm E} \left\lbrace \left\| X_{n-1}^{w} \right\|_{H}^{2} \right\rbrace}}.
\label{A3:94bb}
\end{eqnarray}

\begin{figure}[H]
   \centering
    \includegraphics[width=17cm,height=10cm]{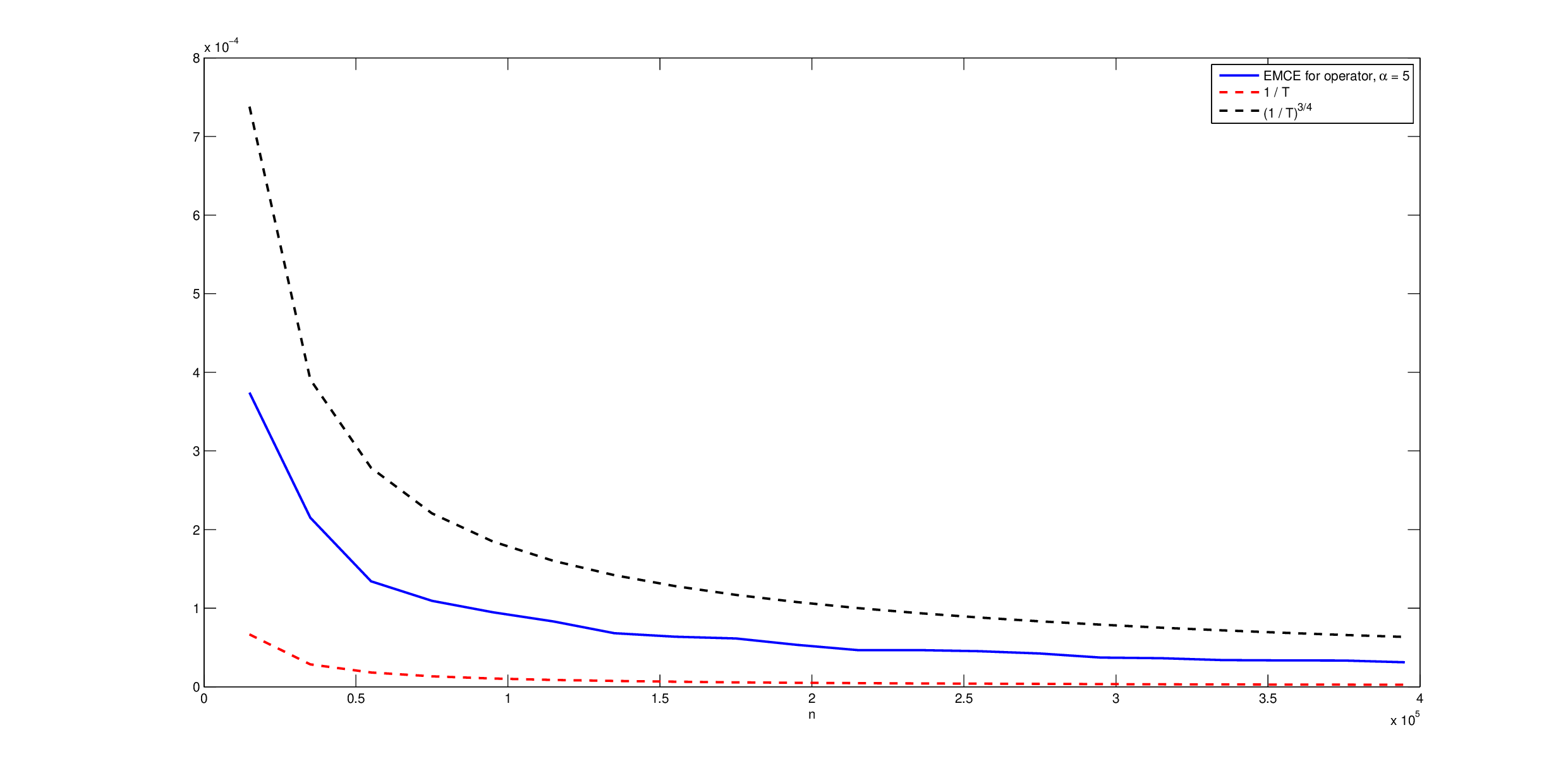}
    \includegraphics[width=17cm,height=10cm]{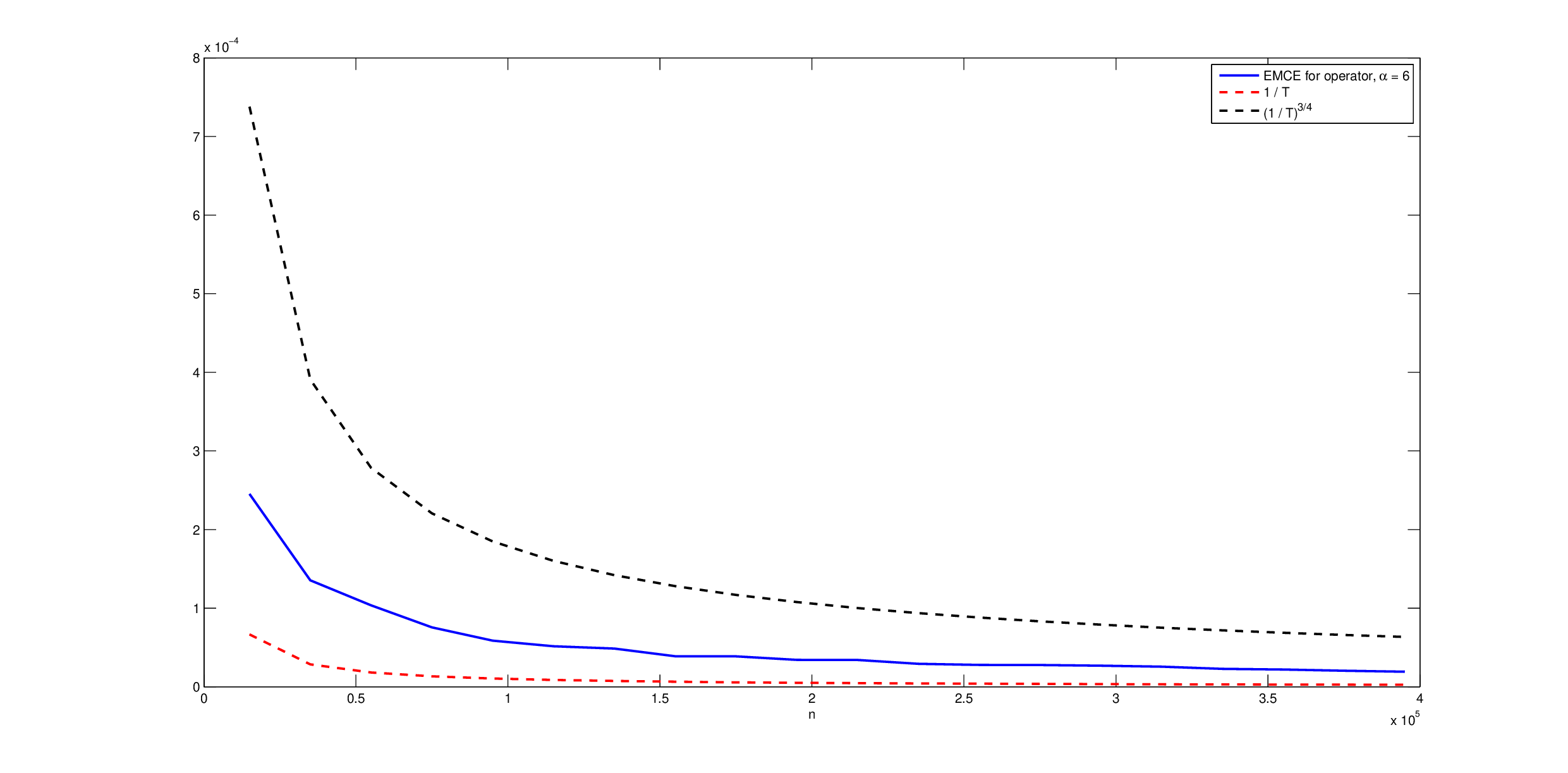}
   
   \vspace{-0.1cm}
\caption[\hspace{0.7cm} Empirical mean square estimation errors of our diagonal approach for large sample sizes and different truncation parameters.]{\small{ ${\rm EMSE}_{\widehat{\rho}_{k_n}}$ values (blue line), in
(\ref{A3:92bb})--(\ref{A3:93bb}), based on $N=700$ simulations, for $\gamma_1 = 0.4$ and $\gamma_2 = 9/20$, considering the sample
sizes   $\left\lbrace n_{t}= 15000+20000(t-1), \ t=1,\dots, 20 \right\rbrace$ and the
corresponding $k_{n,1}$ and $k_{n,2}$ values, for $\alpha_1= 5$
(left-hand side) and $\alpha_2= 6$ (right-hand side), against curves
$(1/n_{t})^{3/4}$ (black dot line) and $1/n_{t}$
(red dot line).}} \label{A3:fig:3}
\end{figure}

\begin{figure}[H]
   \centering
    \includegraphics[width=17cm,height=10cm]{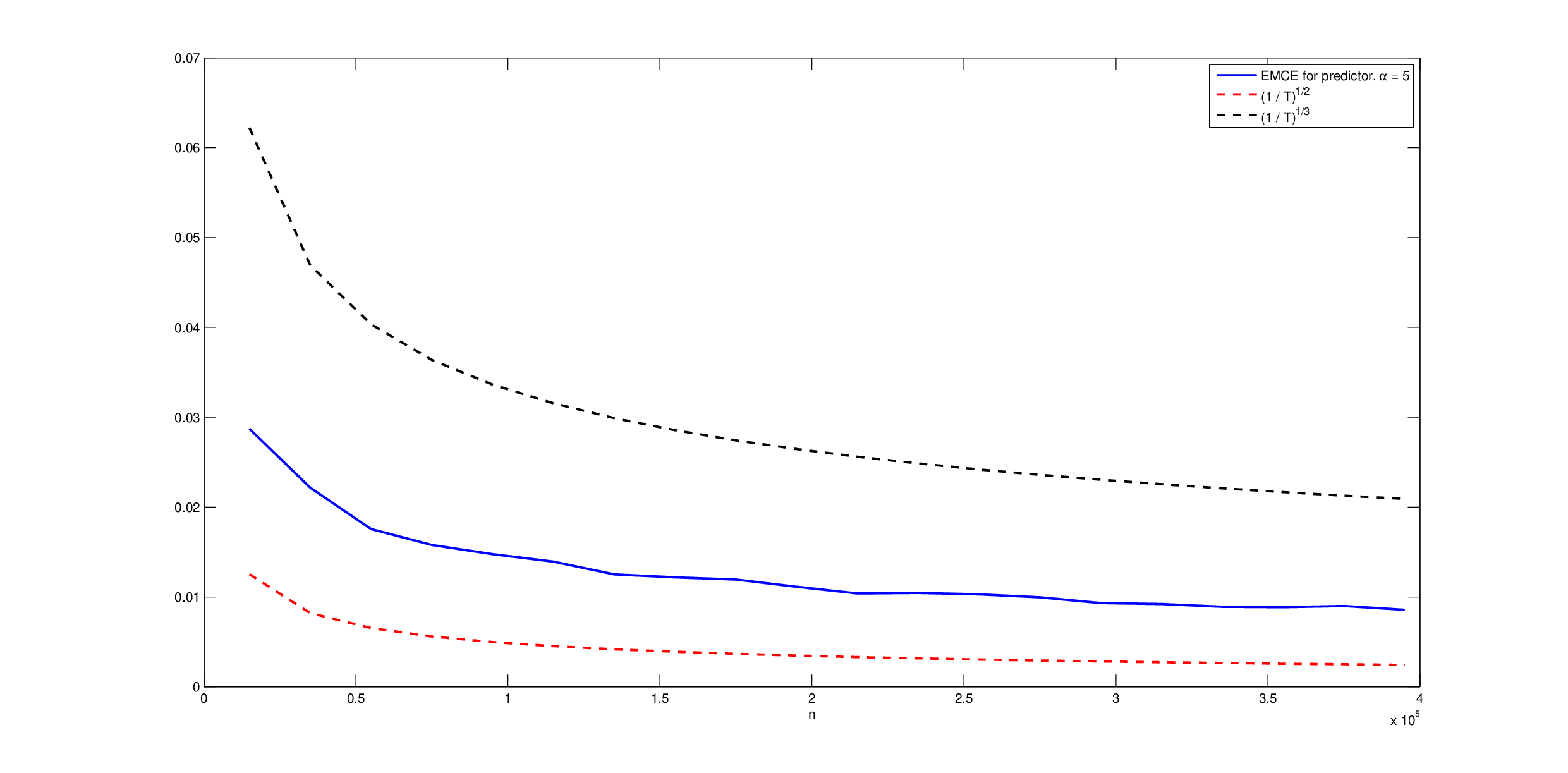}
    \includegraphics[width=17cm,height=10cm]{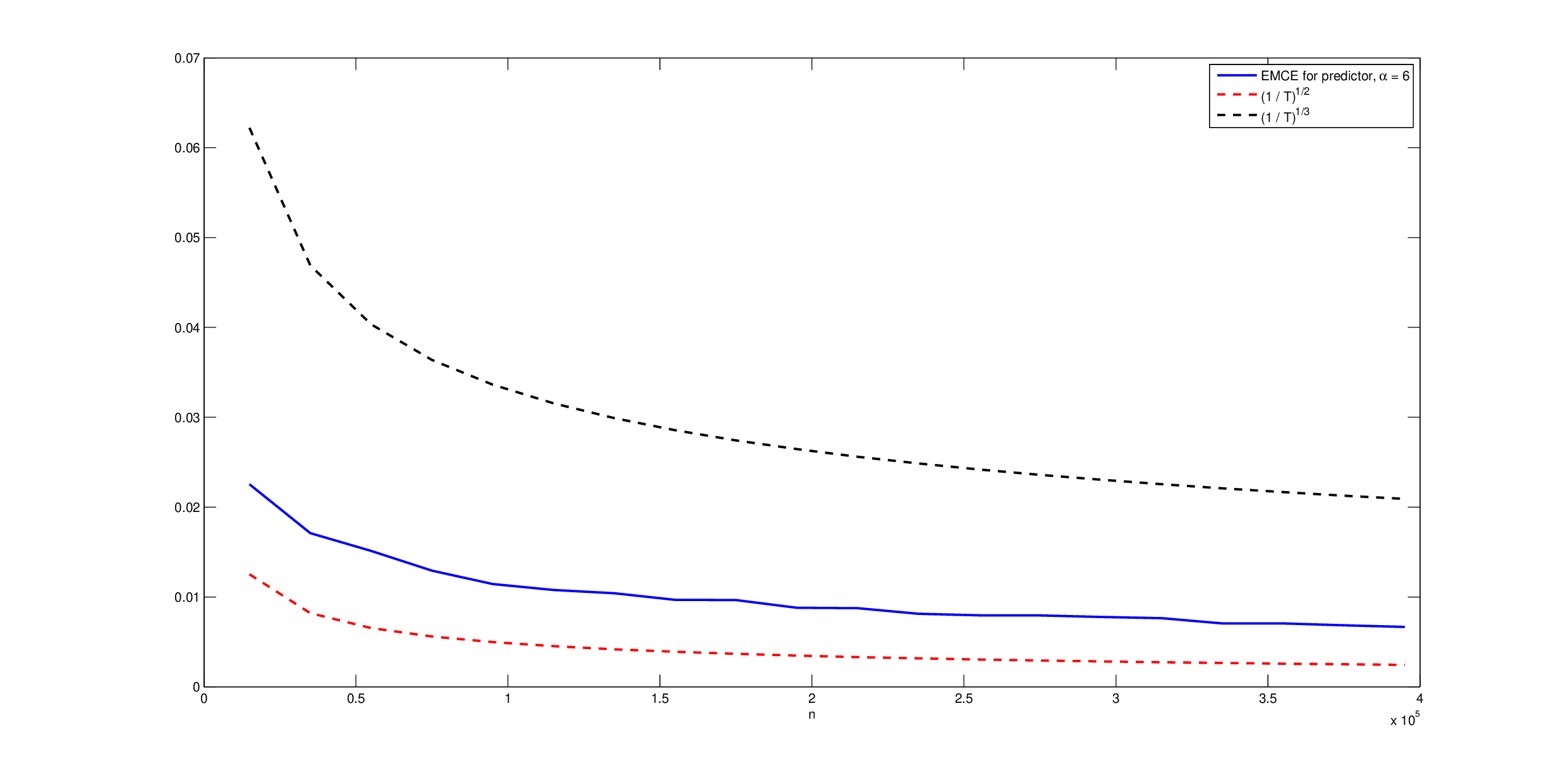}
    
   \vspace{-0.1cm}
\caption[\hspace{0.7cm} Empirical mean square prediction errors of our diagonal approach for large sample sizes and different truncation parameters.]{\small{${\rm UB(EMAE)}_{\widehat{X}_n^{k_n}}$ values (blue line),
in (\ref{A3:94bb}), based on $N=700$ simulations, for $\gamma_1 = 0.4$ and $\gamma_2 = 9/20$, considering the sample
sizes   $\left\lbrace n_{t}= 15000+20000(t-1), \ t=1,\dots, 20 \right\rbrace$ and the
corresponding $k_{n,1}$ and $k_{n,2}$ values, for $\alpha_1= 5$
(left-hand side) and $\alpha_2= 6$ (right-hand side), against curves 
$(1/n_{t})^{1/2}$ (red dot line) and
$(1/n_{t})^{1/3}$ (black dot line).}} \label{A3:fig:4}
\end{figure}

  From  $N = 700$  realizations, for each one of the  elements of the sequence of sample sizes $$\left\lbrace n_t, \ t=1,\dots,20 \right\rbrace = \left\lbrace 15000+20000(t-1), \ t=1,\dots,20 \right\rbrace,$$ the
${\rm EMSE}_{\widehat{\rho}_{k_{n}}}$ and
${\rm UB(EMAE)}_{\widehat{X}_n^{k_n}}$ values, for $\alpha =5$ and
$\alpha =6,$ are displayed in Table \ref{A3:tab:Table1}, where the
abbreviated notations ${\rm MSE}_{\widehat{\rho}_{k_{n,1}}},$ for
${\rm EMSE}_{\widehat{\rho}_{k_{n}}},$ and
${\rm UB}_{\widehat{X}_{n^{k_{n,1}}}},$ for
${\rm UB(EMAE)}_{\widehat{X}_n^{k_n}},$ are used (see also Figures
\ref{A3:fig:3}--\ref{A3:fig:4}).

 \begin{table}[H] 
 \caption[\hspace{0.7cm} Empirical mean square errors of our diagonal approach for large sample sizes and different truncation parameters.]{{\small ${\rm EMSE}_{\widehat{\rho}_{k_n}}$ (here, ${\rm MSE}_{\widehat{\rho}_{k_{n,i}}}$), and
${\rm UB(EMAE)}_{\widehat{X}_n^{k_n}}$ (here,
${\rm UB}_{\widehat{X}_{n^{k_{n,i}}}}$) values,  in
(\ref{A3:92bb})--(\ref{A3:94bb}),  based on $N=700$ simulations, for $\gamma_1 = 0.4$ and $\gamma_2 = 9/20$, considering the
sample sizes   $\left\lbrace n_{t}= 15000+20000(t-1),\ t=1,\dots, 20\right\rbrace$ and the
corresponding $k_{n,1}$ and $k_{n,2}$ values, for $\alpha_1= 5$ and
$\alpha_2= 6$.}}
 \centering
\begin{small}
\begin{tabular}{|c||c|c|c||c|c|c|}
  \hline 
  $n$ & $k_{n,1}$  & ${\rm MSE}_{\widehat{\rho}_{k_{n,1}}}$ & ${\rm UB}_{\widehat{X}_{n^{k_{n,1}}}}$ & $k_{n,2}$  & ${\rm MSE}_{\widehat{\rho}_{k_{n,2}}}$ & ${\rm UB}_{\widehat{X}_{n^{k_{n,2}}}}$\\
  \hline \hline
   $n_1 = 15000$ & $6$ & $3.74 \left( 10\right)^{-4}$ & $2.87 \left( 10\right)^{-2}$ & $4$  & $2.45 \left( 10\right)^{-4}$  & $2.25 \left( 10\right)^{-2}$ \\
  \hline
   $n_2 = 35000$ & $8$ & $2.15 \left( 10\right)^{-4}$ & $2.21 \left( 10\right)^{-2}$ & $5$  & $1.35 \left( 10\right)^{-4}$  & $1.71 \left( 10\right)^{-2}$ \\
  \hline
   $n_3 = 55000$ & $8$ & $1.34 \left( 10\right)^{-4}$ & $1.75 \left( 10\right)^{-2}$ & $6$ & $1.03 \left( 10\right)^{-4}$ & $1.51\left( 10\right)^{-2}$ \\
  \hline
     $n_4 = 75000$ & $9$ & $1.09 \left( 10\right)^{-4}$ & $1.57 \left( 10\right)^{-2}$  & $6$ & $7.55\left( 10\right)^{-5}$ & $1.29 \left( 10\right)^{-2}$ \\
  \hline
     $n_5 = 95000$  & $9$ & $9.48 \left( 10\right)^{-5}$ & $1.47 \left( 10\right)^{-2}$  & $6$ & $5.86 \left( 10\right)^{-5}$ & $1.14 \left( 10\right)^{-2}$ \\
  \hline
     $n_6 = 115000$ & $10$ & $8.31 \left( 10\right)^{-5}$ & $1.39 \left( 10\right)^{-2}$  & $6$ & $5.16 \left( 10\right)^{-5}$ & $1.07 \left( 10\right)^{-2}$ \\
  \hline
    $n_7 = 135000$ & $10$ & $6.81 \left( 10\right)^{-5}$ & $1.25 \left( 10\right)^{-2}$  & $7$ & $4.86 \left( 10\right)^{-5}$ & $1.04 \left( 10\right)^{-2}$ \\
  \hline
      $n_8 = 155000$ & $10$ & $6.37 \left( 10\right)^{-5}$ & $1.21 \left( 10\right)^{-2}$  & $7$ & $3.88\left( 10\right)^{-5}$ & $9.66\left( 10\right)^{-3}$ \\
  \hline
    $n_9 = 175000$ & $11$ & $6.14 \left( 10\right)^{-5}$ & $1.19 \left( 10\right)^{-2}$  & $7$ & $3.87 \left( 10\right)^{-5}$ & $9.65 \left( 10\right)^{-3}$ \\
  \hline
   $n_{10} = 195000$ & $11$ & $5.34 \left( 10\right)^{-5}$ & $1.11 \left( 10\right)^{-2}$  & $7$ & $3.42\left( 10\right)^{-5}$ & $8.79 \left( 10\right)^{-3}$ \\
  \hline
     $n_{11} = 215000$ & $11$ & $4.67 \left(10\right)^{-5}$ & $1.03 \left( 10\right)^{-2}$  & $7$ & $3.40 \left( 10\right)^{-5}$ & $8.74 \left( 10\right)^{-3}$ \\
  \hline
     $n_{12} = 235000$  & $11$ & $4.66 \left( 10\right)^{-5}$ & $1.03 \left( 10\right)^{-2}$  & $7$ & $2.92\left( 10\right)^{-5}$ & $8.12 \left( 10\right)^{-3}$  \\
  \hline
     $n_{13} = 255000$ & $12$ & $4.53 \left( 10\right)^{-5}$ & $1.02\left( 10\right)^{-2}$  & $7$ & $2.77 \left( 10\right)^{-5}$ & $7.95 \left( 10\right)^{-3}$ \\
  \hline
    $n_{14} = 275000$ & $12$ & $4.24 \left(10\right)^{-5}$ & $9.95 \left( 10\right)^{-3}$  & $8$ & $2.77 \left( 10\right)^{-5}$ & $7.94 \left( 10\right)^{-3}$  \\
  \hline
      $n_{15}= 295000$ & $12$ & $3.72 \left( 10\right)^{-5}$ & $9.32 \left( 10\right)^{-3}$ & $8$ & $2.67 \left( 10\right)^{-5}$ & $7.76 \left( 10\right)^{-3}$  \\
  \hline
   $n_{16} = 315000$ & $12$ & $3.62 \left( 10\right)^{-5}$ & $9.21 \left( 10\right)^{-3}$  &  $8$ & $2.55 \left( 10\right)^{-5}$ & $7.64 \left( 10\right)^{-3}$  \\
  \hline
   $n_{17} = 335000$ & $12$ & $3.39 \left(10\right)^{-5}$ & $8.91 \left( 10\right)^{-3}$  & $8$ & $2.28 \left( 10\right)^{-5}$ & $7.04 \left( 10\right)^{-3}$  \\
  \hline
     $n_{18} = 355000$ & $12$ & $3.34 \left( 10\right)^{-5}$ & $8.86 \left( 10\right)^{-3}$ & $8$ & $2.20\left( 10\right)^{-5}$ & $7.04\left( 10\right)^{-3}$  \\
  \hline
     $n_{19} = 375000$  & $13$ & $3.34 \left( 10\right)^{-5}$ & $8.86 \left( 10\right)^{-3}$  & $8$ & $2.04 \left( 10\right)^{-5}$ & $6.84 \left( 10\right)^{-3}$  \\
  \hline
     $n_{20} = 395000$ & $13$ & $3.12\left( 10\right)^{-5}$ & $8.56 \left( 10\right)^{-3}$  & $8$  & $1.92 \left( 10\right)^{-5}$ & $6.65 \left( 10\right)^{-3}$  \\
  \hline
\end{tabular}
\end{small}
  \label{A3:tab:Table1}
\end{table}

\bigskip

In this paper, a one--parameter model of $k_n$ is selected depending on parameter $\alpha$. In \cite[Example 2]{Guillas01}, in
the same spirit, for an equivalent spectral class of operators $C$, a three--parameter model is established for $k_n$ to ensure convergence in quadratic mean in the space $\mathcal{L}(H)$ of the componentwise estimator of $\rho$ constructed from the known
eigenvectors of $C$. The numerical results displayed in Table \ref{A3:tab:Table1} and Figures \ref{A3:fig:3}--\ref{A3:fig:4} illustrate the fact that the proposed componentwise estimator
 $\widehat{\rho}_{k_{n}}$ presents a speed of convergence to $\rho,$ in quadratic mean in $S(H),$  faster than $n^{-1/3},$ which corresponds to the optimal case for the componentwise estimator of $\rho$ proposed in \cite{Guillas01}, in the case of known eigenvectors of $C$; see, in particular, \cite[Theorem 1, Remark 2 and Example 2]{Guillas01}. For larger values of the parameters $\gamma_{1}$ than $2.4,$ and $\alpha $ than $6$,  a faster velocity of convergence of $\widehat{\rho}_{k_{n}}$ to $\rho,$ in quadratic mean in the space $S(H),$ will be obtained. However, larger sample sizes are required for larger values of $\alpha,$ in order to estimate a given number of coefficients of $\rho.$ A more detailed discussion about comparison of the rates of convergence of the ARH(1) plug--in predictors proposed in \cite{AntoniadisSapatinas03,Besseetal00,Bosq00,Guillas01} can be found in the next section.

\textcolor{Crimson}{\subsection{A comparative study}
\label{A3:ACS}}

In this section, the performance of our approach is compared with those ones given in   \cite{AntoniadisSapatinas03,Besseetal00,Bosq00,Guillas01}, including the case of unknown eigenvectors of $C.$  In the last case,  our approach and the approaches presented in  \cite{Bosq00,Guillas01} are implemented in terms of the empirical eigenvectors.

\textcolor{Crimson}{\subsubsection{Theoretical--eigenvector--based componentwise estimators}
\label{A3:cs}}

Let us first compare the performance of our ARH(1) plug--in  predictor, defined in   (\ref{A3:ARHpred}), and the ones formulated in \cite{Bosq00,Guillas01}, in terms of the theoretical eigenvectors $\left\lbrace \phi_j, \ j \geq 1 \right\rbrace$ of $C.$ Note that, in this first part of our comparative study, we consider the previous generated  Gaussian ARH(1) process, with autocovariance and autocorrelation operators defined from equations (\ref{A3:eigvC}) and (\ref{A3:75}), for different rates of convergence to zero of parameters $C_{j}$ and $\rho_{j}^{2},$ $j\geq 1,$ with both sequences being summable sequences.   Since we  restrict our attention to the Gaussian case, conditions A$_{1},$ B$_{1}$ and C$_{1},$ formulated in \cite[pp. 211--212]{Bosq00} are satisfied by the generated ARH(1) process. Similarly, Conditions  H$_{1}$--H$_{3}$ in  \cite[p. 283]{Guillas01}
are satisfied as well.

In \cite[Section 8.2]{Bosq00} the following estimator of $\rho$ is proposed 
\begin{eqnarray}
\widehat{\rho}_n (x) &=& \left(\Pi^{k_n} D_n \widehat{C}_{n}^{-1} \Pi^{k_n}\right)(x) = \displaystyle \sum_{l=1}^{k_n} \widehat{\rho}_{n,l} (x) \phi_l, \quad x \in H, \label{A3:93} \\
\widehat{\rho}_{n,l} (x) &=&  \frac{1}{n-1} \displaystyle \sum_{i=0}^{n-2} \displaystyle \sum_{j=1}^{k_n} \frac{1}{\widehat{C}_{n,j}} \langle \phi_j, x \rangle_{H} X_{i,j}  X_{i+1,l},  \label{A3:94}
\end{eqnarray}
\noindent in the finite dimensional subspace $$H_{k_n} = {\rm span} \left(\phi_1,\dots,\phi_{k_n} \right)$$ of $H,$ where $\Pi^{k_n}$ is the orthogonal projector over $H_{k_n},$ and, as before, $X_{i,j}=\left\langle X_{i},\phi_{j}\right\rangle_{H},$ for $j\geq 1.$ 

A modified estimator of $\rho$ is studied in \cite[Section 2]{Guillas01}, given by
\begin{eqnarray}
\widehat{\rho}_{n,a} (x) &=& \left(\Pi^{k_n} D_n \widehat{C}_{n,a}^{-1} \Pi^{k_n}\right)(x) = \displaystyle \sum_{l=1}^{k_n} \widehat{\rho}_{n,a,l} (x) \phi_l, \quad x \in H, \label{A3:96}\\
\widehat{\rho}_{n,a,l} (x) &=& \frac{1}{n-1} \displaystyle \sum_{i=1}^{n-1} \displaystyle \sum_{j=1}^{k_n} \frac{1}{\displaystyle \max \left(\widehat{C}_{n,j},~a_n \right)} \langle \phi_j, x \rangle_H X_{i,j} X_{i+1,l}\textcolor{green}{,} \label{A3:97}
\end{eqnarray}

\noindent where $$\widehat{C}_{n,a}^{-1}(x) = \displaystyle \sum_{j=1}^{k_n} \frac{1}{\displaystyle \max \left(\widehat{C}_{n,j},~a_n \right)} \langle \phi_j, x \rangle_{H} \phi_j~a.s.$$ Here, $\left\lbrace a_n,\ n \in \mathbb{N} \right\rbrace$ is such that (see  \cite[Theorem 1]{Guillas01})
\begin{equation}
\alpha \frac{C_{k_n}^{\gamma}}{n^{\varepsilon}} \leq a_n \leq \beta \lambda_{k_n}, \quad \alpha >0, \quad 0 < \beta < 1, \quad \varepsilon < 1/2,\quad \gamma \geq 1. \nonumber
\end{equation}

Tables \ref{A3:tab:Table2}--\ref{A3:tab:Table3}  display  the truncated, for two different $k_n$ rules, empirical values of \linebreak${\rm E} \left\lbrace \|\rho \left(X_{n-1} \right)-\widehat{\rho}_{k_{n}}(X_{n-1})\|_{H} \right\rbrace,$ based on $N=700$ generations of each one of the functional samples considered with sizes $n_{t}= 15000+20000(t-1),$ $t=1,\dots, 20,$ when $$C_{j} = b_{C} j^{-\delta_1}, \quad b_{C} > 0, \quad \rho_{j}^{2} = b_{\rho} j^{-\delta_2}, \quad b_{\rho} > 0.$$  Specifically, $\widehat{\rho}_{k_n}$ is computed from equations (\ref{A3:24})--(\ref{A3:25}) (see third column), $\widehat{\rho}_{k_n} = \widehat{\rho}_n$, with $\widehat{\rho}_n$ being given in equations (\ref{A3:93})--(\ref{A3:94}) (see fourth column), and $\widehat{\rho}_{k_n} = \widehat{\rho}_{n,a}$, with $\widehat{\rho}_{n,a}$ being defined in (\ref{A3:96})--(\ref{A3:97})  (see fifth column).

In Table \ref{A3:tab:Table2},  $\delta_1 = 2.4$    $\delta_2 = 1.1,$ and $k_{n} = \lceil n^{1/\alpha} \rceil ,$  for $\alpha= 6,$  according to our \textcolor{Aquamarine}{\textbf{Assumption A3}}, which is also considered in \cite[p. 217]{Bosq00} to ensure weak consistency of the proposed estimator of $\rho$. In Table \ref{A3:tab:Table3}, the same empirical values are displayed for $\delta_{1}=\frac{61}{60},$ $\delta_{2}=1.1,$ and  $k_{n}$ is selected according to  \cite[Example 2]{Guillas01}. Thus, in Table \ref{A3:tab:Table3}, 
\begin{equation}
k_{n}= \lceil n^{\frac{1-2\epsilon}{\delta_{1}(4+2\gamma)}} \rceil,\quad \gamma \geq 1,\ \epsilon<1/2.\label{A3:trGuillas}
\end{equation}

In particular we have chosen $\gamma =2,$ and $\epsilon=0.04\delta_{1}.$  Note that, from  \cite[Theorem 1  and Remark 1]{Guillas01}, for the choice made of  $k_{n}$ in Table \ref{A3:tab:Table3}, convergence to $\rho,$ in quadratic mean in the space $\mathcal{L}(H),$ holds for $\widehat{\rho}_{n,a}$ given in (\ref{A3:96})--(\ref{A3:97}).

\begin{table}[H]
\caption[\hspace{0.7cm} Comparative study on the consistency when eigenvectors are known, for large sample sizes for our truncation parameter.]{\small{Truncated empirical values of ${\rm E} \|\rho \left(X_{n-1}\right)-\widehat{\rho}_{k_{n}}(X_{n-1})\|_{H},$ for $\widehat{\rho}_{k_{n}}$ given in  equations (\ref{A3:24})-(\ref{A3:25}) (third column),  in equations (\ref{A3:93})--(\ref{A3:94}) (fourth column), and  in equations (\ref{A3:96})--(\ref{A3:97}) (fifth column), based on $N=700$ simulations, for $\delta_1 = 2.4$ and $\delta_2 = 1.1,$ considering the
sample sizes  $\left\lbrace n_{t}= 15000+20000(t-1), \ t=1,\dots, 20 \right\rbrace$ and the
corresponding $k_{n} = \lceil n^{1/\alpha} \rceil$ values, for $\alpha= 6$.}}
 \centering
\begin{small}
\begin{tabular}{|c|c||c|c|c|}
  \hline
  $n$ &  $k_{n}$ & Our Approach & Bosq (2000) & Guillas (2001)\\
  \hline \hline
   $n_1 = 15000$ &  $4$  & $2.25 \left( 10\right)^{-2}$ & $2.57 \left( 10\right)^{-2}$ & $2.36 \left( 10\right)^{-2}$  \\
  \hline
   $n_2 = 35000$ &  $5$ & $1.71 \left( 10\right)^{-2}$ & $1.72 \left( 10\right)^{-2}$ & $1.84 \left( 10\right)^{-2}$ \\
  \hline
   $n_3 = 55000$ &  $6$ & $1.51\left( 10\right)^{-2}$ & $1.65\left( 10\right)^{-2}$ & $1.53 \left( 10\right)^{-2}$  \\
  \hline
     $n_4 = 75000$ &  $6$ & $1.29 \left( 10\right)^{-2}$ & $1.46 \left( 10\right)^{-2}$ & $1.37 \left( 10\right)^{-2}$  \\
  \hline
     $n_5 = 95000$  &  $6$ & $1.14 \left( 10\right)^{-2}$ & $1.20 \left( 10\right)^{-2}$ & $1.16 \left( 10\right)^{-2}$  \\
  \hline
     $n_6 = 115000$ & $6$ & $1.07 \left( 10\right)^{-2}$ & $1.10 \left( 10\right)^{-2}$ & $1.11 \left( 10\right)^{-2}$ \\
  \hline
    $n_7 = 135000$ &  $7$ & $1.04 \left( 10\right)^{-2}$ & $1.06 \left( 10\right)^{-2}$ & $1.07 \left( 10\right)^{-2}$ \\
  \hline
      $n_8 = 155000$ &  $7$ & $9.66\left( 10\right)^{-3}$ & $9.91\left( 10\right)^{-3}$ & $1.01 \left( 10\right)^{-2}$ \\
  \hline
    $n_9 = 175000$ &  $7$ & $9.65 \left( 10\right)^{-3}$ & $9.79 \left( 10\right)^{-3}$ & $9.68 \left( 10\right)^{-3}$  \\
  \hline
   $n_{10} = 195000$ &  $7$ & $8.79 \left( 10\right)^{-3}$ & $9.12 \left( 10\right)^{-3}$ & $8.93 \left( 10\right)^{-3}$ \\
  \hline
   $n_{11} = 215000$ &  $7$ & $8.74 \left( 10\right)^{-3}$ & $8.79 \left( 10\right)^{-3}$ & $8.83 \left( 10\right)^{-3}$ \\
  \hline
     $n_{12} = 235000$   & $7$ & $8.12 \left( 10\right)^{-3}$   & $8.69 \left( 10\right)^{-3}$  &$8.75 \left( 10\right)^{-3}$ \\
  \hline
     $n_{13} = 255000$ & $7$ & $7.95 \left( 10\right)^{-3}$ & $8.53 \left( 10\right)^{-3}$ & $8.73 \left( 10\right)^{-3}$ \\
  \hline
    $n_{14} = 275000$  & $8$ & $7.94 \left( 10\right)^{-3}$ & $8.52 \left( 10\right)^{-3}$  & $8.58 \left( 10\right)^{-3}$  \\
  \hline
      $n_{15}= 295000$  & $8$ & $7.76 \left( 10\right)^{-3}$ & $8.49 \left( 10\right)^{-3}$  & $8.36 \left( 10\right)^{-3}$  \\
  \hline
   $n_{16} = 315000$   & $8$  & $7.64 \left( 10\right)^{-3}$ & $7.88 \left( 10\right)^{-3}$  & $8.13 \left( 10\right)^{-3}$ \\
  \hline
   $n_{17} = 335000$ & $8$ & $7.04 \left( 10\right)^{-3}$ & $7.24 \left( 10\right)^{-3}$  & $7.59 \left( 10\right)^{-3}$ \\
  \hline
     $n_{18} = 355000$ & $8$ & $7.04\left( 10\right)^{-3}$ & $7.23\left( 10\right)^{-3}$  & $6.92 \left( 10\right)^{-3}$ \\
  \hline
     $n_{19} = 375000$  & $8$ & $6.84 \left( 10\right)^{-3}$ & $6.89 \left( 10\right)^{-3}$  & $6.90 \left( 10\right)^{-3}$ \\
  \hline
     $n_{20} = 395000$ & $8$ & $6.65 \left( 10\right)^{-3}$   & $6.67 \left( 10\right)^{-3}$ & $6.85 \left( 10\right)^{-3}$ \\
  \hline
\end{tabular}
\end{small}
  \label{A3:tab:Table2}
\end{table}

\begin{table}[H] 
\caption[\hspace{0.7cm} Comparative study on the consistency when eigenvectors are known, for large sample sizes for truncation parameter in Guillas (2001).]{\small{Truncated empirical values of ${\rm E} \|\rho \left(X_{n-1}\right) -\widehat{\rho}_{k_{n}}(X_{n-1})\|_{H},$ for $\widehat{\rho}_{k_{n}}$ given in  equations (\ref{A3:24})--(\ref{A3:25}) (third column),  in equations (\ref{A3:93})--(\ref{A3:94}) (fourth column), and  in equations (\ref{A3:96})--(\ref{A3:97}) (fifth column), based on $N=700$ simulations, for $\delta_1 = \frac{61}{60}$ and $\delta_2 = 1.1,$ considering the
sample sizes  $\left\lbrace n_{t}= 15000+20000(t-1), \ t=1,\dots, 20 \right\rbrace$ and the
corresponding $k_{n}$ given in (\ref{A3:trGuillas}).}}
\centering
\begin{small}
\begin{tabular}{|c|c||c|c|c|}
  \hline 
  $n$ & $k_{n}$   & Our Approach & Bosq (2000) & Guillas (2001)\\
  \hline \hline
   $n_1 = 15000$ & $2$   & $9.91\left( 10\right)^{-3}$ & $1.39 \left( 10\right)^{-2}$ & $1.26 \left( 10\right)^{-2}$ \\
  \hline
   $n_2 = 35000$ & $3$  & $8.78\left( 10\right)^{-3}$& $1.34 \left( 10\right)^{-2}$    & $1.24 \left( 10\right)^{-2}$ \\
  \hline
   $n_3 = 55000$ & $3$   & $7.89\left( 10\right)^{-3}$ & $1.15 \left( 10\right)^{-2}$   & $1.14 \left( 10\right)^{-2}$ \\
  \hline
     $n_4 = 75000$ & $3$    & $6.49\left( 10\right)^{-3}$ & $1.01 \left( 10\right)^{-2}$  & $8.58 \left( 10\right)^{-3}$ \\
  \hline
     $n_5 = 95000$  & $3$ & $6.36\left( 10\right)^{-3}$  & $9.09 \left( 10\right)^{-3}$    & $8.29 \left( 10\right)^{-3}$ \\
  \hline
     $n_6 = 115000$ &$ 3 $& $6.14\left( 10\right)^{-3}$ & $7.65 \left( 10\right)^{-3}$    & $7.26\left( 10\right)^{-3}$  \\
  \hline
    $n_7 = 135000$ & $3 $ & $5.91\left( 10\right)^{-3}$ & $7.03 \left( 10\right)^{-3}$  & $6.69\left( 10\right)^{-3}$ \\
  \hline
      $n_8 = 155000$ &$ 3 $ & $5.73\left( 10\right)^{-3}$ & $6.77 \left( 10\right)^{-3}$   & $6.54\left( 10\right)^{-3}$  \\
  \hline
    $n_9 = 175000$ & $3 $ & $5.44\left( 10\right)^{-3}$& $6.74 \left( 10\right)^{-3}$  & $6.16\left( 10\right)^{-3}$   \\
  \hline
   $n_{10} = 195000$ & $3$   & $5.10\left( 10\right)^{-3}$   & $6.69 \left( 10\right)^{-3}$  & $5.97\left( 10\right)^{-3}$ \\
  \hline
   $n_{11} = 215000$ & $4$ & $5.01\left( 10\right)^{-3}$ & $6.48 \left( 10\right)^{-3}$    & $5.94\left( 10\right)^{-3}$  \\
  \hline
     $n_{12} = 235000$  & $4 $   &$4.85\left( 10\right)^{-3}$ & $6.45 \left( 10\right)^{-3}$  & $5.83\left( 10\right)^{-3}$ \\
  \hline
     $n_{13} = 255000$ & $4 $   & $4.17\left( 10\right)^{-3}$ & $6.17\left( 10\right)^{-3}$  & $5.68\left( 10\right)^{-3}$  \\
  \hline
    $n_{14} = 275000$ & $4$   & $4.64\left( 10\right)^{-3}$ & $5.99 \left( 10\right)^{-3}$  & $5.60\left( 10\right)^{-3}$ \\
  \hline
      $n_{15}= 295000$ & $4$    & $4.55\left( 10\right)^{-3}$  & $5.94 \left( 10\right)^{-3}$ & $5.58\left( 10\right)^{-3}$  \\
  \hline
   $n_{16} = 315000$ & $4$  & $4.48\left( 10\right)^{-3}$ & $5.69 \left( 10\right)^{-3}$  & $5.50\left( 10\right)^{-3}$   \\
  \hline
   $n_{17} = 335000$ & $4 $ & $4.38\left( 10\right)^{-3}$ & $5.58 \left( 10\right)^{-3}$    & $5.44\left( 10\right)^{-3}$  \\
  \hline
     $n_{18} = 355000$ & $4 $  & $4.16\left( 10\right)^{-3}$  & $5.45 \left( 10\right)^{-3}$& $5.42\left( 10\right)^{-3}$  \\
  \hline
     $n_{19} = 375000$  & $4$  & $3.91\left( 10\right)^{-3}$ & $5.34\left( 10\right)^{-3}$   & $5.32\left( 10\right)^{-3}$  \\
  \hline
     $n_{20} = 395000$ & $4 $  & $3.86\left( 10\right)^{-3}$ & $5.29 \left( 10\right)^{-3}$ & $5.26\left( 10\right)^{-3}$ \\
  \hline
\end{tabular}
\end{small}
  \label{A3:tab:Table3}
\end{table}

One can observe in Table \ref{A3:tab:Table2} a similar performance of the three methods compared with the truncation order kn satisfying
\textcolor{Aquamarine}{\textbf{Assumption A3}}, with slightly worse results being obtained from the estimator defined in (\ref{A3:96})--(\ref{A3:97}), specially, for the
sample size $n_{8}=155000.$
Furthermore,  in Table \ref{A3:tab:Table3}, a better performance of our approach is observed for the smallest sample
sizes (from $n_{1}=15000$ until $n_{4} = 75000$). For the remaining largest sample sizes, only slight differences are observed, with,
again, a better performance of our approach, very close to the other two approaches presented in \cite{Bosq00,Guillas01}.

\textcolor{Crimson}{\subsubsection{Empirical--eigenvector--based componentwise estimators}}

In this section, we address the case where $\left\lbrace \phi_{j},\ j\geq 1 \right\rbrace$ are unknown, as is often the case in practice. Specifically, for a given sample size $n$, let $\left\lbrace \phi_{n,j}, \ j \geq 1 \right\rbrace$  be the empirical counterpart of the theoretical eigenvectors $\left\lbrace \phi_{j},\ j\geq 1 \right\rbrace$, satisfying, for every $j \geq 1$,
\begin{equation}
C_n \left( \phi_{n,j} \right) = \frac{1}{n} \displaystyle \sum_{i=0}^{n-1} \langle X_i, \phi_{n,j} \rangle_H X_i = C_{n,j} \phi_{n,j}, \nonumber
\end{equation}

\noindent where $\left\lbrace C_{n,j}, \ j \geq 1 \right\rbrace$  denotes the system of eigenvalues associated with the system of empirical eigenvectors \linebreak $\left\lbrace \phi_{n,j}, \ j \geq 1 \right\rbrace$. We then consider the following estimators for comparison purposes

\begin{eqnarray}
\widetilde{\rho}_{n,j} &=& \frac{\frac{1}{n-1} \displaystyle \sum_{i=0}^{n-2} \widetilde{X}_{i,j} \widetilde{X}_{i+1,j}}{\frac{1}{n} \displaystyle \sum_{i=0}^{n-1} \left( \widetilde{X}_{i,j} \right)^2},\quad \widetilde{\rho}_{k_{n}}=\sum_{j=1}^{k_{n}}\widetilde{\rho}_{n,j}\phi_{n,j}\otimes \phi_{n,j},\label{A3:e1rho}\\
\widetilde{\rho}_n (x) &=& \left(\widetilde{\Pi}^{k_n} D_n C_{n}^{-1} \widetilde{\Pi}^{k_n} \right)(x) = \displaystyle \sum_{l=1}^{k_n} \widetilde{\rho}_{n,l} (x) \phi_{n,l}, \quad x \in H,\nonumber \\
\widetilde{\rho}_{n,l} (x) &=&  \frac{1}{n-1} \displaystyle \sum_{i=0}^{n-2} \displaystyle \sum_{j=1}^{k_n} \frac{1}{C_{n,j}} \langle \phi_{n,j}, x \rangle_{H} \widetilde{X}_{i,j}  \widetilde{X}_{i+1,l},\label{A3:eqrhobosq}\\
\widetilde{\rho}_{n,a} (x) &=& \left(\widetilde{\Pi}^{k_n} D_n C_{n,a}^{-1} \widetilde{\Pi}^{k_n} \right)(x) = \displaystyle \sum_{l=1}^{k_n} \widetilde{\rho}_{n,a,l} (x) \phi_{n,l}, \quad x \in H,\nonumber \\
\widetilde{\rho}_{n,a,l} (x) &=&  \frac{1}{n-1} \displaystyle \sum_{i=0}^{n-2} \displaystyle \sum_{j=1}^{k_n} \frac{1}{\displaystyle \max \left(C_{n,j},a_n \right)} \langle \phi_{n,j}, x \rangle_{H} \widetilde{X}_{i,j}  \widetilde{X}_{i+1,l,}\label{A3:eqrGuillas}
\end{eqnarray}
\noindent where, for $i\in \mathbb{Z},$ and $j\geq 1,$ $\widetilde{X}_{i,j}=\left\langle X_{i}, \phi_{n,j}\right\rangle_{H},$ $\widetilde{\Pi}^{k_n}$ denotes the orthogonal projector into the space $$\widetilde{H}_{k_n} = {\rm span} \left(\phi_{n,1}, \dots, \phi_{n,k_n} \right).$$

The Gaussian ARH(1) process is generated under \textcolor{Aquamarine}{\textbf{Assumptions A1--A2}}, as well as $C_{1}^{\prime }$   in \cite[p. 218]{Bosq00}. Note that  conditions $A_{1}$ and $B_{1}^{\prime}$ in \cite{Bosq00} already hold. Moreover, as given in \cite[Theorem 8.8 and Example 8.6]{Bosq00}, for $$C_j = b_{C} j^{-\delta_1}, \quad b_{C} > 0, \quad \delta_1 > 0,$$ \noindent with, in particular, $\delta_1 = 2.4,$ and for $$\rho_j = b_{\rho} j^{-\delta_2}, \quad b_{\rho} > 0,$$ with $\delta_2 = 1.1,$, the estimator $\widetilde{\rho}_n$ converges almost surely to $\rho$ under the condition $$\frac{nC_{kn}^{2}}{\ln(n)\left(\displaystyle \sum_{j=1}^{k_{n}}b_{j}\right)^{2}}\longrightarrow \infty,$$ where $$b_1 = 2 \sqrt{2} \left(C_1 - C_2 \right)^{-1}, \quad b_{j} = 2\sqrt{2}\max\left\lbrace (C_{j-1}-C_{j})^{-1}, (C_{j}-C_{j+1})^{-1}\right\rbrace,~j\geq 2.$$ In Table \ref{A3:tab:Table4}, $k_{n} = \lceil \ln(n) \rceil$ has been tested; see \cite[Example 8.6]{Bosq00}.

\bigskip

\begin{table}[H]
\caption[\hspace{0.7cm} Comparative study on the consistency when eigenvectors are unknown and truncation in Bosq (2000) is used.]{\small{Truncated empirical values of ${ \rm E} \left\lbrace \left\| \rho \left(X_{n-1} \right) - \widetilde{\rho}_{k_n} \left(X_{n-1} \right) \right\|_{H} \right\rbrace$, for $\widetilde{\rho}_{k_n} =  \widetilde{\rho}_{k_n}$ given in equation (\ref{A3:e1rho}) (third column), $\widetilde{\rho}_{k_n} = \widetilde{\rho}_{n}$ defined in equation (\ref{A3:eqrhobosq}) (fourth column) and $\widetilde{\rho}_{k_n} = \widetilde{\rho}_{n,a}$ defined in equation (\ref{A3:eqrGuillas}) (fifth column), based on $N=700$ simulations, for $\delta_1 = 2.4$ and $\delta_2 = 1.1,$ considering the sample sizes $\left\lbrace n_{t}= 15000+20000(t-1), \ t=1,\dots, 20 \right\rbrace$ and  $k_{n} = \lceil \ln(n) \rceil$.}}
\centering
\begin{small}
\begin{tabular}{|c|c||c|c|c|}
  \hline
  $n$ & $k_{n}$  & Our approach  & Bosq (2000) & Guillas (2001)\\
  \hline \hline
   $n_1 = 15000$ & $9$ & $8.42 \left( 10\right)^{-2}$ & $1.061$ & $1.035$ \\
  \hline
   $n_2 = 35000$ & $10$  & $5.51 \left( 10\right)^{-2}$ & $1.019$& $1.005$  \\
  \hline
   $n_3 = 55000$ & $10$   & $4.75 \left( 10\right)^{-2}$ & $1.017$& $0.999$ \\
  \hline
     $n_4 = 75000$ & $11 $& $4.43 \left( 10\right)^{-2}$ & $1.015$ & $0.995$ \\
  \hline
     $n_5 = 95000$  & $11$ & $3.68 \left( 10\right)^{-2}$ & $1.013$ & $0.988$ \\
  \hline
     $n_6 = 115000$ & $11$  & $3.51 \left( 10\right)^{-2}$ & $1.011$ & $0.963$ \\
  \hline
    $n_7 = 135000$ &$ 11 $ & $3.23 \left( 10\right)^{-2}$ & $1.008$ & $0.925$  \\
  \hline
      $n_8 = 155000$ & $11 $& $2.95 \left( 10\right)^{-2}$  & $1.007$  & $0.912$  \\
  \hline
    $n_9 = 175000$ &$ 12$  & $2.94 \left( 10\right)^{-2}$ & $1.006$  & $0.911$  \\
  \hline
   $n_{10} = 195000$ & $12$  & $2.80 \left( 10\right)^{-2}$ & $0.995$ & $0.891$  \\
  \hline
   $n_{11} = 215000$ & $12$ & $2.71 \left( 10\right)^{-2}$ & $0.902$  & $0.862$  \\
  \hline
     $n_{12} = 235000$  &$ 12$  & $2.59 \left( 10\right)^{-2}$  & $0.890$  & $0.820$  \\
  \hline
     $n_{13} = 255000$ & $12$ & $2.58 \left( 10\right)^{-2}$ & $0.878$  &  $0.800$  \\
  \hline
    $n_{14} = 275000$ & $12$ & $2.35 \left( 10\right)^{-2}$ & $0.872$   & $0.783$ \\
  \hline
      $n_{15}= 295000$ & $12 $ &  $2.28 \left( 10\right)^{-2}$ & $0.860$   & $0.778$  \\
  \hline
   $n_{16} = 315000$ & $12$  & $2.27 \left( 10\right)^{-2}$ & $0.842$ & $0.747$  \\
  \hline
   $n_{17} = 335000$ & $12 $ & $2.16\left( 10\right)^{-2}$ & $0.822$ & $0.714$  \\
  \hline
     $n_{18} = 355000$ & $12$ & $2.14\left( 10\right)^{-2}$  & $0.800$ & $0.707$ \\
  \hline
     $n_{19} = 375000$  & $12$  & $2.09\left( 10\right)^{-2}$  & $0.778$& $0.687$  \\
  \hline
     $n_{20} = 395000$ & $12$ & $2.06\left( 10\right)^{-2}$ & $0.769$ & $0.662$  \\
  \hline
\end{tabular}
\end{small}
  \label{A3:tab:Table4}
\end{table}

\bigskip

  A better performance of our estimator (\ref{A3:e1rho}) in comparison with estimator (\ref{A3:eqrhobosq}), formulated in \cite{Bosq00},   and estimator (\ref{A3:eqrGuillas}), formulated in \cite[Example 4 and Remark 4]{Guillas01}, is observed in Table \ref{A3:tab:Table4}. Note that, in particular,  in \cite[Example 4 and Remark 4]{Guillas01}, smaller values of $k_{n}$ than $\ln(n)$ are required for a given sample size $n,$ to ensure convergence in quadratic mean, and, in particular, weak--consistency.  However, considering a smaller discretization step size $\Delta t = 0.015$ than in Table \ref{A3:tab:Table4}, where $\Delta t = 0.08$, and for $k_{n}= \lceil n^{1/6} \rceil,$ (i.e., $\alpha =6$), we obtain in  Table \ref{A3:tab:Table5},   for the same parameter values $\delta_1 = 2.4$ and $\delta_2 = 1.1,$ better results than in Table \ref{A3:tab:Table4}, since a smaller number of coefficients of $\rho$ (parameters) to be estimated is considered in Table \ref{A3:tab:Table5}, from a  richer sample information (coming from the  smaller discretization step size considered). One can also observe in Table \ref{A3:tab:Table5} a similar performance of the three approaches studied. In Table \ref{A3:tab:Table6}, the value $k_{n} = \lceil e' n^{1/\left(8 \delta_1 + 2 \right)} \rceil$, with $e' = \frac{17}{10}$  proposed in \cite[Example 4 and Remark 4]{Guillas01}  is considered to compute the truncated empirical values of ${ \rm E} \left\lbrace  \|\rho(X_{n-1})-\widetilde{\rho}_{k_{n}}(X_{n-1})\|_{H} \right\rbrace,$ for $\widetilde{\rho}_{k_{n}}$ defined in equation  (\ref{A3:e1rho}) (third column), for \linebreak $\widetilde{\rho}_{k_{n}}= \widetilde{\rho}_n$ given in  equation (\ref{A3:eqrhobosq})  (fourth column), and for $\widetilde{\rho}_{k_{n}}=
\widetilde{\rho}_{n,a}$ in equation (\ref{A3:eqrGuillas}) (fifth column). A similar performance of the three
approaches is observed, with the exception of $n_{20} = 395000,$ where the approach presented in \cite{Guillas01} displays a slightly better
performance

\bigskip

  \begin{table}[H] 
  \caption[\hspace{0.7cm} Comparative study on the consistency when eigenvectors are unknown and our truncation parameter is used, with a small discretization step.]{\small{Truncated empirical values of ${ \rm E} \left\lbrace \left\| \rho \left(X_{n-1} \right) - \widetilde{\rho}_{k_n} \left(X_{n-1} \right) \right\|_{H} \right \rbrace$, for $\widetilde{\rho}_{k_{n}}$ defined in equation  (\ref{A3:e1rho}) (third column), for $\widetilde{\rho}_{k_{n}}= \widetilde{\rho}_n$ given in  equation (\ref{A3:eqrhobosq})  (fourth column), and for $\widetilde{\rho}_{k_{n}}=
\widetilde{\rho}_{n,a}$ in equation (\ref{A3:eqrGuillas})
(fifth column), based on $N=200$  (due to high-dimensionality) simulations, for $\delta_1 = 2.4$ and $\delta_2 = 1.1,$ considering the sample sizes $\left\lbrace n_{t}= 15000+20000(t-1),\ t=1,\dots, 20 \right\rbrace$ and  $k_{n}= \lceil n^{1/6} \rceil.$}}
  \centering
\begin{small}
\begin{tabular}{|c|c||c|c|c|}
  \hline
  $n$ & $k_{n}$  &  Our approach  & Bosq (2000) & Guillas (2001)\\
  \hline \hline
   $n_1 = 15000$ & $4$ & $9.88 \left( 10\right)^{-2}$ & $9.25\left( 10\right)^{-2}$ & $0.106$ \\
  \hline
   $n_2 = 35000$ & $5$  & $9.52 \left( 10\right)^{-2}$ & $9.07 \left( 10\right)^{-2}$ & $9.86 \left( 10\right)^{-2}$ \\
  \hline
   $n_3 = 55000$ & $6 $  & $9.12\left( 10\right)^{-2}$ & $8.92 \left( 10\right)^{-2}$ & $9.39  \left( 10\right)^{-2}$\\
  \hline
     $n_4 = 75000$ & $6$ & $8.48 \left( 10\right)^{-2}$ & $8.64 \left( 10\right)^{-2}$& $8.98 \left( 10\right)^{-2}$\\
  \hline
     $n_5 = 95000$  & $6 $& $7.61 \left( 10\right)^{-2}$ & $8.30 \left( 10\right)^{-2}$& $8.46 \left( 10\right)^{-2}$ \\
  \hline
     $n_6 = 115000$ & $6 $& $7.05 \left( 10\right)^{-2}$  & $7.96\left( 10\right)^{-2}$ & $8.04 \left( 10\right)^{-2}$ \\
  \hline
    $n_7 = 135000$ &$ 7$ & $6.99 \left( 10\right)^{-2}$  & $7.84\left( 10\right)^{-2}$ & $7.82\left( 10\right)^{-2}$ \\
  \hline
      $n_8 = 155000$ & $7$ & $6.70 \left( 10\right)^{-2}$  & $7.45\left( 10\right)^{-2}$ & $7.40\left( 10\right)^{-2}$ \\
  \hline
    $n_9 = 175000$ &$ 7$ & $6.49 \left( 10\right)^{-2}$  & $7.03\left( 10\right)^{-2}$ & $7.07\left( 10\right)^{-2}$ \\
  \hline
   $n_{10} = 195000$ & $7$  & $5.88 \left( 10\right)^{-2}$  & $6.74\left( 10\right)^{-2}$ & $6.80\left( 10\right)^{-2}$ \\
  \hline
   $n_{11} = 215000$ & $7$ & $5.63 \left( 10\right)^{-2}$  & $6.46\left( 10\right)^{-2}$  & $6.57\left( 10\right)^{-2}$ \\
  \hline
     $n_{12} = 235000$  &$ 7 $ & $5.30 \left( 10\right)^{-2}$  & $6.28\left( 10\right)^{-2}$  &  $6.37\left( 10\right)^{-2}$ \\
  \hline
     $n_{13} = 255000$ & $7$ & $5.05 \left( 10\right)^{-2}$  & $6.19\left( 10\right)^{-2}$  &  $6.24\left( 10\right)^{-2}$ \\
  \hline
    $n_{14} = 275000$ &$ 8$ & $4.88 \left( 10\right)^{-2}$  & $5.99\left( 10\right)^{-2}$     &  $6.15\left( 10\right)^{-2}$ \\
  \hline
      $n_{15}= 295000$ & $8 $ & $4.58 \left( 10\right)^{-2}$  & $5.74\left( 10\right)^{-2}$  &  $6.04\left( 10\right)^{-2}$ \\
  \hline
   $n_{16} = 315000$ & $8$  & $4.24 \left( 10\right)^{-2}$  & $5.52\left( 10\right)^{-2}$ & $5.93\left( 10\right)^{-2}$ \\
  \hline
   $n_{17} = 335000$ & $8$  & $3.86 \left( 10\right)^{-2}$  & $5.24\left( 10\right)^{-2}$ & $5.70\left( 10\right)^{-2}$ \\
  \hline
     $n_{18} = 355000$ &$ 8$ & $3.70 \left( 10\right)^{-2}$ & $5.02\left( 10\right)^{-2}$  & $5.53\left( 10\right)^{-2}$ \\
  \hline
     $n_{19} = 375000$  &$ 8 $ & $3.55 \left( 10\right)^{-2}$  & $4.88\left( 10\right)^{-2}$ & $5.36\left( 10\right)^{-2}$ \\
  \hline
     $n_{20} = 395000$ &$ 8$ & $3.46 \left( 10\right)^{-2}$ & $4.70\left( 10\right)^{-2}$  & $5.23 \left( 10\right)^{-2}$ \\
  \hline
\end{tabular}
\end{small}
  \label{A3:tab:Table5}
\end{table}

  \begin{table}[H]
\caption[\hspace{0.7cm} Comparative study on the consistency when eigenvectors are unknown and truncation in Guillas (2001) is used.] {\small{Truncated empirical values of ${ \rm E} \left\lbrace \left\| \rho \left(X_{n-1} \right) - \widetilde{\rho}_{k_n} \left(X_{n-1} \right) \right\|_{H} \right\rbrace$, for $\widetilde{\rho}_{k_{n}}$ defined in equation  (\ref{A3:e1rho}) (third column), for $\widetilde{\rho}_{k_{n}}= \widetilde{\rho}_n$ given in  equation (\ref{A3:eqrhobosq})  (fourth column), and for $\widetilde{\rho}_{k_{n}}=
\widetilde{\rho}_{n,a}$ in equation (\ref{A3:eqrGuillas})
(fifth column), based on $N=200$  (due to high-dimensionality) simulations, for $\delta_1 = 2.4$ and $\delta_2 = 1.1,$ considering the sample sizes $\left\lbrace n_{t}= 15000+20000(t-1), \ t=1,\dots, 20 \right\rbrace$ and  $k_{n} = \lceil e' n^{1/\left(8 \delta_1 + 2 \right)} \rceil,~e' = \frac{17}{10}.$}}
   \centering
\begin{small}
\begin{tabular}{|c|c||c|c|c|}
  \hline
  $n$ & $k_{n}$  &  Our approach  & Bosq (2000) & Guillas (2001)  \\
  \hline \hline
   $n_1 = 15000$ & $2$ & $6.78 \left( 10\right)^{-2}$ & $8.77\left( 10\right)^{-2}$ & $6.64 \left( 10\right)^{-2}$\\
  \hline
   $n_2 = 35000$ & $2$ & $6.72 \left( 10\right)^{-2}$ & $8.61\left( 10\right)^{-2}$  & $6.30 \left( 10\right)^{-2}$ \\
  \hline
   $n_3 = 55000$ & $2$ & $6.46 \left( 10\right)^{-2}$ & $8.48\left( 10\right)^{-2}$  & $6.17\left( 10\right)^{-2}$ \\
  \hline
     $n_4 = 75000$ & $2$ & $6.24 \left( 10\right)^{-2}$ & $8.20\left( 10\right)^{-2}$ & $5.76\left( 10\right)^{-2}$ \\
  \hline
     $n_5 = 95000$  & $2$ & $5.42 \left( 10\right)^{-2}$ &  $7.84\left( 10\right)^{-2}$ & $5.03\left( 10\right)^{-2}$ \\
  \hline
     $n_6 = 115000$ & $2$ & $4.84 \left( 10\right)^{-2}$ & $7.34\left( 10\right)^{-2}$ & $4.56\left( 10\right)^{-2}$\\
  \hline
    $n_7 = 135000$ & $2$ & $4.27\left( 10\right)^{-2}$ & $6.95\left( 10\right)^{-2}$ & $3.94\left( 10\right)^{-2}$ \\
  \hline
      $n_8 = 155000$ & $2 $& $3.64\left( 10\right)^{-2}$ & $6.60\left( 10\right)^{-2}$ & $3.65\left( 10\right)^{-2}$ \\
  \hline
    $n_9 = 175000$ &$ 3$ & $3.51\left( 10\right)^{-2}$ & $6.52\left( 10\right)^{-2}$ & $3.42\left( 10\right)^{-2}$ \\
  \hline
   $n_{10} = 195000$ & $3 $ & $3.38\left( 10\right)^{-2}$ & $6.16\left( 10\right)^{-2}$ & $3.24\left( 10\right)^{-2}$ \\
  \hline
   $n_{11} = 215000$ & $3$ & $3.16\left( 10\right)^{-2}$ & $5.78\left( 10\right)^{-2}$ & $2.85\left( 10\right)^{-2}$ \\
  \hline
     $n_{12} = 235000$  &$ 3$  & $2.98\left( 10\right)^{-2}$ & $5.53\left( 10\right)^{-2}$ & $2.60\left( 10\right)^{-2}$ \\
  \hline
     $n_{13} = 255000$ & $3$ & $2.83\left( 10\right)^{-2}$ & $5.15\left( 10\right)^{-2}$ & $2.34\left( 10\right)^{-2}$  \\
  \hline
    $n_{14} = 275000$ &$ 3$ & $2.50 \left( 10\right)^{-2}$ & $4.85\left( 10\right)^{-2}$  & $2.05 \left( 10\right)^{-2}$ \\
  \hline
      $n_{15}= 295000$ & $3$  & $2.23\left( 10\right)^{-2}$ & $4.46\left( 10\right)^{-2}$ & $1.83 \left( 10\right)^{-2}$  \\
  \hline
   $n_{16} = 315000$ & $3$ & $2.15\left( 10\right)^{-2}$ & $4.30\left( 10\right)^{-2}$ & $1.58\left( 10\right)^{-2}$  \\
  \hline
   $n_{17} = 335000$ & $3$  & $2.06 \left( 10\right)^{-2}$ & $4.14\left( 10\right)^{-2}$ & $1.40\left( 10\right)^{-2}$\\
  \hline
     $n_{18} = 355000$ & $3$ & $1.98 \left( 10\right)^{-2}$ & $3.95 \left( 10\right)^{-2}$ & $1.24 \left( 10\right)^{-2}$ \\
  \hline
     $n_{19} = 375000$  & $3$  & $1.89 \left( 10\right)^{-2}$ & $3.77 \left( 10\right)^{-2}$ & $1.05 \left( 10\right)^{-2}$ \\
  \hline
     $n_{20} = 395000$ & $ 3$ & $1.82 \left( 10\right)^{-2}$ & $3.70 \left( 10\right)^{-2}$ & $9.93 \left( 10\right)^{-3}$ \\
  \hline
\end{tabular}
\end{small}
  \label{A3:tab:Table6}
\end{table}

  \textcolor{Crimson}{\subsubsection{Kernel--based nonparametric and penalized estimation}}

In practice, curves are observed in discrete times, and should be approximated by smooth functions. In  \cite{Besseetal00},  the following optimization problem
is considered:
\begin{equation}
\widehat{X}_{i} = argmin \left\| L \widehat{X}_{i} \right\|_{L^2}^{2},~\widehat{X}_{i}(t_j) = X_i (t_j),\quad j=1,\dots,p,~i=0,\dots,n-1,
\label{A3:sdata}
\end{equation}
\noindent where $L$ is a linear differential operator of order $d.$ Our interpolation is computed by Matlab \emph{smoothingspline} method.  Non-linear kernel regression is then considered, in terms of the smoothed functional data, solution to (\ref{A3:sdata}), as follows:
\begin{eqnarray}
\widehat{X}_{n}^{h_n} &=& \widehat{\rho}_{h_n}(X_{n-1}), \quad
\widehat{\rho}_{h_n} (x) =\frac{\displaystyle \sum_{i=0}^{n-2} \widehat{X}_{i+1}K\left(\frac{ \left\| \widehat{X}_{i} - x \right\|_{L^2}^{2}}{h_n} \right)}{\displaystyle \sum_{i=0}^{n-2} K\left(\frac{ \left\| \widehat{X}_{i} - x \right\|_{L^2}^{2}}{h_n} \right)}, \nonumber
\end{eqnarray}

\noindent where $K$ is the usual Gaussian kernel, and $$\left\| \widehat{X}_{i} - x \right\|_{L^2}^{2}=\int (\widehat{X}_{i}(t)-x(t))^{2}dt,\quad i=0,\dots,n-2.$$

Alternatively, in \cite{Besseetal00}, prediction, in the context of functional autoregressive processes (FAR(1) processes), under the linear assumption on $\rho,$ which is considered to be a compact operator, with  $\|\rho\|<1,$ is also studied, from smooth data $\widehat{X}_{1},\dots, \widehat{X}_{n},$ solving the optimization problem

 \begin{equation}
\displaystyle \min_{\widehat{X}_i \in H_q} \frac{1}{n} \displaystyle \sum_{i=0}^{n-1} \left( \frac{1}{p} \displaystyle \sum_{j=1}^{p} \left(X_i (t_j) -  \widehat{X}_{i}^{q,l} (t_j)\right)^2 + l \left\| D^2 \widehat{X}_{i}^{q,l} \right\|_{L^2}^{2} \right),\label{A3:esr}
\end{equation}

\noindent where $l$ is the smoothing parameter, $H_q$ is the $q$--dimensional functional subspace spanned by the leading eigenvectors of the autocovariance operator $C$ associated with its largest eigenvalues. Thus, smoothness and rank constraint are considered in the computation of the solution to the optimization problem (\ref{A3:esr}). Such a solution is obtained by means of functional PCA. 

The following regularized empirical estimators of $C$ and $D$ are then considered, with inversion of $C$
in the subspace $H_q$:
$$\widehat{C}_{q,l}=\frac{1}{n}\sum_{i=0}^{n-1}\widehat{X}_{i}\otimes \widehat{X}_{i},\quad \widehat{D}_{q,l}=\frac{1}{n-1}\sum_{i=0}^{n-2}\widehat{X}_{i}\otimes \widehat{X}_{i+1}.$$

Thus, the regularized estimator of $\rho$ is given by
$$\widehat{\rho}_{q,l}=\widehat{D}_{q,l}\widehat{C}_{q,l}^{-1},$$
 and the predictor $$\widehat{X}_{n}^{q,l}=\widehat{\rho}_{q,l}X_{n-1}.$$
Due to computational cost limitations, in Table  \ref{A3:tab:Table7}, the following statistics are evaluated  to compare the performance of the two above-referred  prediction methodologies:
\begin{equation}
EMAE_{\widehat{X}_n}^{h_n} = \frac{1}{p} \displaystyle \sum_{j=1}^{p} \left(X_{n}(t_j) - \widehat{X}_{n}^{h_n}(t_j) \right)^2, \label{A3:109}
\end{equation}
\begin{equation}
EMAE_{\widehat{X}_n}^{q,l} = \frac{1}{p} \displaystyle \sum_{j=1}^{p} \left(X_{n}(t_j) - \widehat{X}_{n}^{q,l}(t_j) \right)^2. \label{A3:113}
\end{equation}

 \begin{table}[H] 
 \caption[\hspace{0.7cm} Comparative study on the consistency when Besse et al. (2000) proposal is tested, for small sample sizes.]{\small{$EMAE_{\widehat{X}_n}^{h_{n,i}},$ $i=1,2,$  and $EMAE_{\widehat{X}_n}^{q,l}$ values (see (\ref{A3:109}) and (\ref{A3:113}), respectively), with $q = 7$,
based on $N=200$ simulations, for $\delta_1 = 2.4$ and $\delta_2 = 1.1,$ considering now the sample sizes  $\left\lbrace n_{t}= 750+500(t-1), \ t=1,\dots, 13 \right\rbrace$  $h_{n,1} = 0.1$ and $h_{n,2} = 0.3$.}}
 \centering
\begin{small}
\begin{tabular}{|c||c|c||c|}
  \hline
  $n$   &  $EMAE_{\widehat{X}_n}^{h_{n,1}}$ & $EMAE_{\widehat{X}_n}^{h_{n,2}}$ & $EMAE_{\widehat{X}_n}^{q,l}$ \\
  \hline \hline
   $n_1 = 750$ &  $8.57 \left(10 \right)^{-2}$ & $8.85 \left(10 \right)^{-2}$ & $8.99 \left(10 \right)^{-2}$\\
  \hline
   $n_2 = 1250$  & $7.67 \left(10 \right)^{-2}$ & $8.43 \left(10 \right)^{-2}$ & $8.69 \left(10 \right)^{-2}$\\
  \hline
   $n_3 = 1750$   & $7.15 \left(10 \right)^{-2}$ & $7.12 \left(10 \right)^{-2}$ & $8.05 \left(10 \right)^{-2}$\\
  \hline
     $n_4 = 2250$   & $7.09 \left(10 \right)^{-2}$& $6.87 \left(10 \right)^{-2}$ & $7.59 \left(10 \right)^{-2}$\\
  \hline
     $n_5 = 2750$   & $6.87 \left(10 \right)^{-2}$ & $6.67 \left(10 \right)^{-2}$ & $7.31 \left(10 \right)^{-2}$\\
  \hline
     $n_6 = 3250$   & $6.52 \left(10 \right)^{-2}$& $5.92 \left(10 \right)^{-2}$ & $7.28 \left(10 \right)^{-2}$\\
  \hline
    $n_7 = 3750$   & $6.20 \left(10 \right)^{-2}$& $5.56 \left(10 \right)^{-2}$ & $7.13 \left(10 \right)^{-2}$\\
  \hline
      $n_8 = 4250$    & $6.06 \left(10 \right)^{-2}$ & $5.32 \left(10 \right)^{-2}$ & $7.06 \left(10 \right)^{-2}$\\
  \hline
    $n_9 = 4750$  & $5.67 \left(10 \right)^{-2}$ & $5.25 \left(10 \right)^{-2}$ & $6.47 \left(10 \right)^{-2}$\\
  \hline
   $n_{10} = 5250$    & $5.24 \left(10 \right)^{-2}$ & $5.12 \left(10 \right)^{-2}$ & $6.08 \left(10 \right)^{-2}$\\
  \hline
   $n_{11} = 5750$    & $5.01 \left(10 \right)^{-2}$ & $4.82 \left(10 \right)^{-2}$& $5.75 \left(10 \right)^{-2}$\\
  \hline
     $n_{12} = 6250$    & $4.90 \left(10 \right)^{-2}$ &$4.49 \left(10 \right)^{-2}$ & $5.33 \left(10 \right)^{-2}$\\
  \hline
     $n_{13} = 6750$ & $4.87 \left(10 \right)^{-2}$ & $3.87 \left(10 \right)^{-2}$ & $4.97 \left(10 \right)^{-2}$\\
  \hline
\end{tabular}
\end{small}
  \label{A3:tab:Table7}
\end{table}

\bigskip

It can be observed  a similar performance of the kernel--based and penalized FAR(1) predictors, from smooth functional data, which is also comparable, considering one realization, to the performance obtained in Table  \ref{A3:tab:Table6}, from the empirical eigenvectors.

\textcolor{Crimson}{\subsubsection{Wavelet--based prediction for ARH(1) processes}}

The approach presented in \cite{AntoniadisSapatinas03} is now studied.  Specifically, wavelet-based regularization is applied to obtain smooth estimates of
the sample paths. The projection onto the space $V_{J},$ generated by translations of the scaling function $\phi_{Jk},\ k=0,\dots, 2^{J}-1,$ at level $J,$ associated with a multiresolution analysis of $H,$ is first considered. For a given primary resolution level $j_{0}$, with $j_{0}<J,$ the following wavelet decomposition at $J-j_{0}$ resolution levels can be computed for any projected curve $\Phi_{V_{J}}X_{i},$ in the space $V_{J},$ for  $i=0,\dots,n-1:$
\begin{eqnarray}
& & \Phi_{V_{J}}X_{i} = \displaystyle \sum_{k=0}^{2^{j_0}-1} c_{j_0 k}^{i} \phi_{j_0 k} + \displaystyle \sum_{j=j_0}^{J-1}  \displaystyle \sum_{k=0}^{2^j - 1} d_{jk}^{i} \psi_{jk},\nonumber\\ & &  c_{j_0 k}^{i} = \langle \Phi_{V_{J}}X_i, \phi_{j_0 k} \rangle_H,~d_{j k}^{i} = \langle \Phi_{V_{J}}X_i, \psi_{j k} \rangle_H. \nonumber 
\end{eqnarray}

For $i=0,\dots,n-1,$ the following variational problem is solved to obtain the  smooth estimate of
the curve  $X_{i}:$
\begin{equation}
\displaystyle \inf_{f^{i} \in H} \left\lbrace \left\| \Phi_{V_{J}} X_i - f^{i} \right\|_{L^2}^{2} + \lambda \left\| \Phi_{V_{j_0}^{\bot}} f \right\|^2 ;~f \in H \right\rbrace, \label{A3:130}
\end{equation}
\noindent where $\Phi_{V_{j_0}^{\bot}}$ denotes the orthogonal projection operator of $H$ onto the orhogonal complement of $V_{j_0},$ and for $i=0,1\dots n-1,$
 $$f^{i}=\sum_{k=0}^{2^{j_{0}}-1}\alpha_{j_0 k}^{i} \phi_{j_0 k} + \displaystyle \sum_{j=j_0}^{\infty }  \displaystyle \sum_{k=0}^{2^j - 1} \beta _{jk}^{i} \psi_{jk}.$$
 
 Using the equivalent sequence of norms of fractional Sobolev spaces of order $s$ with $s>1/2,$ on a suitable interval (in our case, $s=\delta_1$), the minimization of (\ref{A3:130}) is equivalent to the optimization problem, for $i=0,\dots,n-1,$ \begin{equation}
 \sum_{k=0}^{2^{j_{0}}-1}(\alpha_{j_0 k}^{i}-c_{j_0 k}^{i})^{2}+\displaystyle \sum_{j=j_0}^{J-1}  \displaystyle \sum_{k=0}^{2^j - 1}(d_{jk}^{i}-\beta _{jk}^{i} )^{2}+ \sum_{j=j_0}^{\infty }  \displaystyle \sum_{k=0}^{2^j - 1}\lambda 2^{js}[\beta _{jk}^{i}]^{2}.\label{A3:rewav}
 \end{equation}
 The solution to (\ref{A3:rewav}) is given by, for $i=0,\dots,n-1,$
\begin{eqnarray}
\widehat{\alpha_{j_0 k}^{i}} &=& c_{j_0 k}^{i},\quad k=0,1,\dots,2^{j_0}-1, \nonumber \\ 
\widehat{\beta_{j_0 k}^{i}} &=& \frac{d_{jk}^{i}}{(1+\lambda 2^{2sj})},\quad j=j_0,\dots,J-1,~k=0,1,\dots,2^{j}-1, \nonumber \\
\widehat{\beta_{j_0 k}^{i}} &=& 0,\quad j \geq J,~k=0,1,\dots,2^{j}-1. \nonumber 
\end{eqnarray}

 In particular, in the subsequent computations, we have considered the following value of the smoothing parameter $\lambda $ (see \cite{Angelinietal03}):
 $$\widehat{\lambda}^{M} = \frac{\left(\displaystyle \sum_{j=1}^{M} \sigma_{j}^{2} \right) \left(\displaystyle \sum_{j=1}^{M} C_{j}\right)}{n}.$$

 The following smoothed data are then computed
 \begin{equation}
\widetilde{X}_{i,\widehat{\lambda}^{M}} = \displaystyle \sum_{k=0}^{2^{j_0}-1} \widehat{\alpha_{j_0 k}^{i}}\phi_{j_0 k} + \displaystyle \sum_{j=j_0}^{J-1}  \displaystyle \sum_{k=0}^{2^j - 1} \widehat{\beta_{j_0 k}^{i}} \psi_{jk}, \nonumber
\end{equation}
\noindent removing the trend $$\widetilde{a}_{n,\widehat{\lambda}^{M}} = \frac{1}{n} \displaystyle \sum_{i=0}^{n-1} \widetilde{X}_{i,\widehat{\lambda}^{M}}$$ to obtain $$\widetilde{Y}_{i,\widehat{\lambda}^{M}}=\widetilde{X}_{i,\widehat{\lambda}^{M}}-\widetilde{a}_{n,\widehat{\lambda}^{M}},\quad i=0,\dots,n-1,$$ for the computation of
\begin{eqnarray}
\widetilde{\rho}_{n,\widehat{\lambda}^{M}} (x) &=& \left(\widetilde{\Pi}_{\widehat{\lambda}^{M}}^{k_n} \widetilde{D}_{n,\widehat{\lambda}^{M}} \widetilde{C}_{n,\widehat{\lambda}^{M}}^{-1} \widetilde{\Pi}_{\widehat{\lambda}^{M}}^{k_n} \right)(x) = \displaystyle \sum_{l=1}^{k_n} \widetilde{\rho}_{n,\widehat{\lambda}^{M},l} (x) \widetilde{\phi}_{l}^{M},~ x \in H,\nonumber\\
\widetilde{\rho}_{n,\widehat{\lambda}^{M},l} (x) &=&   \displaystyle \sum_{j=1}^{k_n}\frac{1}{n-1} \displaystyle \sum_{i=0}^{n-2} \frac{1}{\widetilde{C}_{n,\widehat{\lambda}^{M},j}} \langle \widetilde{\phi}_{j}^{M}, x \rangle_{H} \widetilde{Y}_{i,\widehat{\lambda}^{M},j}  \widetilde{Y}_{i+1,\widehat{\lambda}^{M},l}, \nonumber
\end{eqnarray}
\noindent for $x \in H$ and

\begin{eqnarray}
\widetilde{C}_{n,\widehat{\lambda}^{M}}&=& \frac{1}{n} \displaystyle \sum_{i=0}^{n-1} \widetilde{Y}_{i,\widehat{\lambda}^{M}} \otimes \widetilde{Y}_{i,\widehat{\lambda}^{M}} , \nonumber
\end{eqnarray}

\noindent where $$\widetilde{Y}_{i,\widehat{\lambda}^{M},j}=\left\langle \widetilde{Y}_{i,\widehat{\lambda}^{M}},\widehat{\phi}_{j,\widehat{\lambda}^{M}}\right\rangle,$$ and $$\widetilde{C}_{n,\widehat{\lambda}^{M},j}=\left\langle\widetilde{C}_{n,\widehat{\lambda}^{M}}\widehat{\phi}_{j,
\widehat{\lambda}^{M}}\right\rangle,$$ for every $j \geq 1$. Table \ref{A3:tab:Table8} displays  the empirical truncated approximation of the expectation \linebreak  ${ \rm E} \left\lbrace \|\widetilde{\rho}_{k_{n}} (X_{n-1})-\rho(X_{n-1})\|_{H} \right\rbrace$ and ${ \rm E} \left\lbrace \|\widetilde{\rho}_{n,\widehat{\lambda}^{M}} (X_{n-1})-\rho(X_{n-1})\|_{H} \right\rbrace,$
 respectively obtained applying our approach, and the approach in \cite{AntoniadisSapatinas03}, in the estimation of the autocorrelation operator $\rho$. Here, we have tested $k_{n_{i}}= \lceil n^{1/ \alpha_i} \rceil ,$ $i=1,2,$ with $\alpha_1 = 6,$ according to \textcolor{Aquamarine}{\textbf{Assumption A3}}, and $\alpha_{2}>4\delta_{1},$  according to
 $$H_{4}:~nC_{k_{n}}^{4}\rightarrow \infty$$ in \cite[p. 149]{AntoniadisSapatinas03}. In particular, we have considered   $\delta_1 = 2.4,$ and $\alpha_2 = 10.$ From the results displayed in Table
\ref{A3:tab:Table8}, one can observe a similar performance for the two truncation rules implemented,
and approaches compared, for the small sample sizes tested. A similar accuracy is also displayed by the approaches
presented in \cite{Besseetal00}, for such small sample sizes (see Table \ref{A3:tab:Table7}).

\bigskip

\begin{table}[H] 
\caption[\hspace{0.7cm} Comparative study on the consistency when Antoniadis and Sapatinas (2003) proposal is tested, for small sample sizes.]{\small{Truncated empirical values of ${ \rm E} \left\lbrace \|\rho(X_{n-1})-\widetilde{\rho}_{k_{n}}(X_{n-1})\|_{H} \right\rbrace,$  with $\widetilde{\rho}_{k_{n}}$ defined in equation  (\ref{A3:e1rho}),
and of ${ \rm E} \left\lbrace \|\widetilde{\rho}_{n,\widehat{\lambda}^{M}} (X_{n-1})-\rho(X_{n-1})\|_{H} \right\rbrace,$, based on $N=200$
simulations, for $\delta_1 = 2.4$ and $\delta_2 = 1.1,$ considering  the sample sizes  $\left\lbrace n_{t}= 750+500(t-1), \ t=1,\dots, 13 \right\rbrace$, using $\widehat{\lambda}_M,~M = 50$,
and the corresponding $k_{n,i} = \lceil n^{1/\alpha_i} \rceil,$ for $\alpha_1 = 6$ and $\alpha_2 = 10$. Here, O.A. means \emph{Our Approach}.}}
\centering
\begin{small}
\begin{tabular}{|c|c|c|c||c|c|c|}
  \hline
  $n$   & $k_{n,1}$ & O.A. &  \cite{AntoniadisSapatinas03} & $k_{n,2}$ & O.A. &  \cite{AntoniadisSapatinas03} \\
  \hline \hline
   $ 750$ & 3 & $0.070$ & $0.091$ & 1 & $0.064$ & $0.059$ \\
  \hline
   $ 1250$ & 3 & $0.055$ & $0.087$ & 2 & $0.051$ & $ 0.043$ \\
  \hline
   $1750$ & 3 & $0.047$ & $0.080$ & 2 & $0.045$ & $ 0.039$ \\
  \hline
     $ 2250$ & 3 & $0.041$ & $0.079$ & 2 & $0.041$ & $0.038$ \\
  \hline
     $2750$ & 3 & $0.037$ & $0.073$ & 2 & $0.036$ & $0.035$ \\
  \hline
     $ 3250$ & 3 & $0.034$ & $0.072$ & 2 & $0.033$ & $0.031$\\
  \hline
    $ 3750$  & 3 & $0.033$ & $0.068$ & 2 & $0.033$ & $0.029$ \\
  \hline
      $4250$  & 4 & $0.033$ & $0.067$ & 2 & $0.031$ & $0.029$ \\
  \hline
    $ 4750$ & 4 & $0.032$ & $0.066$ & 2 & $0.031$ & $0.026$ \\
  \hline
   $ 5250$  & 4 & $0.031$ & $0.064$ & 2 & $0.028$ & $0.023$ \\
  \hline
   $5750$  & 4 & $0.030$ & $0.060$ & 2 & $0.020$ & $0.019$ \\
  \hline
     $ 6250$  & 4 & $0.028$ & $0.058$ & 2 & $0.017$ & $0.015$ \\
  \hline
     $ 6750$ & 4 & $0.028$ & $0.056$ & 2 & $0.019$ & $0.014$ \\
  \hline
\end{tabular}
\end{small}
  \label{A3:tab:Table8}
\end{table}

\textcolor{Crimson}{\section{Final comments}
\label{A3:sec:5}}

As noted before, in this paper, the eigenvectors of  $C$ are considered to be known in the derivation of the results on
consistency. This assumption is satisfied, e.g., when the random initial condition is given as the solution, in the mean-square
sense, of a stochastic differential equation driven by white noise (e.g., the Wiener measure), since the eigenvectors of the
differential operator involved in that equation coincide with the eigenvectors of the autocovariance operator of the ARH(1)
process.
In the case where the eigenvectors of the autocovariance operator are unknown, the numerical results displayed in Tables \ref{A3:tab:Table4}--\ref{A3:tab:Table6} illustrate the fact that our approach displays, in terms of the empirical eigenvectors, very
similar prediction results to those obtained with the implementation of the componentwise estimators proposed in \cite{Bosq00,Guillas01},
with a better performance of our approach in the more unfavorable case, corresponding to a large discretization step size,
and truncation order (see Table \ref{A3:tab:Table4} computed for $k_n = \lceil \ln(n) \rceil$).

Regarding \textcolor{Aquamarine}{\textbf{Assumption A2}}, \textcolor{Crimson}{Remark} \ref{A3:rem1def} provides an example where this assumption is satisfied. However, our approach
can still be applied in a wider range of situations. Wavelet bases are well suited for sparse representation of functions;
recent work has considered combining them with sparsity-inducing penalties, both for semiparametric regression (see,
e.g., \cite{WandOrmerod11}), and for regression with functional or kernel predictors (see \cite{WandOrmerod11,Zhaoetal12,Zhaoetal15}, among others). The latter papers focused on $\ell_1$ penalization, also known as the lasso (see \cite{Tibshirani96}), in the wavelet domain. Alternatives to the lasso include the SCAD penalty  by \cite{FanLi01}, and the adaptive lasso by \cite{Zou06}. The $\ell_1$ penalty in the elastic net criterion has the effect of shrinking small coefficients to zero.
This can be interpreted as imposing a prior that favors a sparse estimate.
The above mentioned smoothing techniques, based on wavelets, can be applied to obtain a smooth sparse approximation
 $\widehat{X}_{1},\dots \widehat{X}_{n}$ of the functional data $X_{1},\dots, X_{n},$  whose empirical auto-covariance operator $$\widehat{C}_{n}= \frac{1}{n}\sum_{i=0}^{n-1}\widehat{X}_{i}\otimes \widehat{X}_{i}$$
 and cross-covariance operator $$ \widehat{D}_{n}= \frac{1}{n-1}\sum_{i=0}^{n-2}\widehat{X}_{i}\otimes \widehat{X}_{i+1} $$
 admits a diagonal representation in terms of wavelets.

 In the literature, shrinkage approaches for estimating a high--dimensional covariance matrix are employed
to circumvent the limitations of the sample covariance matrix. In particular, a new family of nonparametric Stein--type
shrinkage covariance estimators is proposed in \cite{Touloumis15} (see also references therein), whose members are written as a convex
linear combination of the sample covariance matrix and of a predefined invertible diagonal target matrix. These results
can be applied to our framework, considering the shrinkage estimators of the autocovariance and cross-covariance operators, with respect to a common suitable wavelet basis, which can lead to an empirical diagonal representation of both
operators.

\bigskip

In the Supplementary Material provided (see \textcolor{Crimson}{Appendix} \ref{A3:Supp}), a numerical example is provided to illustrate the performance of our
approach, in the case of a pseudo--diagonal autocorrelation operator.

\textcolor{Crimson}{\section{Supplementary Material: non--diagonal autocorrelation operator} \label{A3:Supp}}

This Section provides as a numerical example where the methodology proposed in such paper still works beyond the considered \textcolor{Aquamarine}{\textbf{Assumption A2}}. In particular, this section illustrates the performance of the proposed estimation methodology, when \textcolor{Aquamarine}{\textbf{Assumption A2}} is not satisfied, but $\rho$ is close to be diagonal
in some sense. The numerical results obtained are compared to those ones
derived from the computation of the ARH(1) predictors, based on the componentwise estimators proposed in  \cite{Bosq00,Guillas01} where this diagonal assumption is
not required. The Gaussian ARH(1) process generated has autocorrelation
operator $\rho$ with coefficients $\rho_{j,h}$ with respect to  the basis $\left\lbrace \phi_{j}\otimes \phi_{h}, \ j,h\geq 1 \right\rbrace,$ given by
  \begin{eqnarray}
\rho_{j,j}^{2} &=& \left(\frac{\lambda_j \left( \left(- \Delta \right)_{\left(a,b\right)} \right)}{\lambda_1 \left( \left(- \Delta \right)_{\left(a,b\right)} \right) - \epsilon}\right)^{- \delta_2}, \label{A3:118}\end{eqnarray}
\noindent in the diagonal, and outside of the diagonal

 \begin{eqnarray}
\rho_{j,j+a}^{2} &=& \frac{0.01}{5 a^2},~ a=1,\ldots,5, \quad \rho_{j+a,j}^{2} = \frac{0.02}{5 a^2},~ a=1,\ldots,5, \label{A3:119}
\end{eqnarray}

\noindent where  $\rho_{j,j+a}^{2} = \rho_{j+a,j}^{2} = 0$ when $a \geq 6$. The coefficients of the autocovariance operator $C_{\varepsilon}$  of the innovation process $\varepsilon,$ with respect to   the mentioned basis $\left\lbrace \phi_{j}\otimes \phi_{h}, \ j,h\geq 1 \right\rbrace,$ are given by $$\sigma_{j,j}^{2} = C_j \left(1 - \rho_{j,j}^{2} \right)$$ in the diagonal, and outside
of the diagonal by
\begin{eqnarray}
\sigma_{j,j+a}^{2} &=& \frac{0.015}{5 a^2},~ a=1,2,3,4,5, \quad \sigma_{j+a,j}^{2} =\frac{0.01}{5 a^2},~ a=1,2,3,4,5, \label{A3:121}
\end{eqnarray}
\noindent where  $\sigma_{j,j+a}^{2} = \sigma_{j+a,j}^{2} = 0$ when $a \geq 6$. Table \ref{A3:tab:Table11} below displays the empirical truncated values of \linebreak ${ \rm E} \left\lbrace \left\| \rho(X_{n-1}) - \widehat{\rho}_{k_n}^{ND} (X_{n-1}) \right\|_{H}\right\rbrace$ based on $N = 200$
simulations of each one of the 20 functional samples considered, with sizes $\left\lbrace n_t = 15000 + 20000(t - 1), \ t = 1, \ldots, 20 \right\rbrace$, for the corresponding $k_n$ values obtained, in this case, by the rule $k_n = \lceil n^{1/\alpha} \rceil$, with $\alpha = 6$. We have considered parameter $\delta_1 = 2.4$ in the definition of the eigenvalues
of $C$; but, in this case, as noted before, operators $\rho$ and $C_{\varepsilon}$ are non-diagonal
(see equations \ref{A3:119}--\ref{A3:121}). The estimators of $\rho$ and the associated plug--in
predictors are computed, for the three approaches compared, under the assumption that the eigenvectors of C are known.

As expected, in Table \ref{A3:tab:Table11}, an outperformance of the approaches in \cite{Bosq00,Guillas01} is observed in comparison with our methodology. However, for large sample
sizes, the ARH(1) prediction methodology proposed here still can be applied
with an order of magnitude of $10^{-2}$ for the empirical errors associated with $\widehat{\rho}_{k_n}$ given in equation \ref{A3:e1rho}. Thus, in the pseudodiagonal autocorrelation
operator case, in some sense, our approach could still be considered. As referred in our paper, an example is given in the case where the autocovariance
and autocorrelation operators admit a sparse representation in terms of a
suitable orthonormal wavelet basis (see, for instance, \cite{Angelinietal03,AntoniadisSapatinas03}).

\bigskip

\begin{table}[H]
\caption[\hspace{0.7cm} Comparative study on the  consistency when a pseudodiagonal scenario is regarded and eigenvectors are known.]{\small{Truncated empirical values of ${ \rm E} \left\lbrace \left\| \rho(X_{n-1}) - \widehat{\rho}_{k_n}^{ND} (X_{n-1}) \right\|_{H} \right\rbrace$, for $\widehat{\rho}_{k_n}^{ND}$ given in equations  (\ref{A3:24})--(\ref{A3:25}) (third column),  in equations (\ref{A3:93})--(\ref{A3:94}) (fourth column), and  in equations (\ref{A3:96})--(\ref{A3:97}) (fifth column), from the non--diagonal data generated by equations (\ref{A3:118})--(\ref{A3:121}), based on $N=200$ (due to high--dimensionality) simulations, for $\delta_1 = 2.4$ and $\delta_2 = 1.1,$ considering the
sample sizes  $\left\lbrace n_{t}= 15000+20000(t-1), \ t=1,\dots, 20 \right\rbrace$ and the corresponding $k_{n} = \lceil n^{1/\alpha} \rceil,~\alpha = 6$ values. The eigenvectors $\left\lbrace \phi_j, \ j \geq 1 \right\rbrace$ are assumed to be known.}} 
\centering
\begin{small}
\begin{tabular}{|c|c||c|c|c|}
  \hline
  $n$ & $k_{n}$  &    Our approach  & Bosq (2000) & Guillas (2001)\\
  \hline \hline
   $n_1 = 15000$ & $4$ & $0.581$ & $8.94 \left(10 \right)^{-2}$ & $0.1055$\\
  \hline
   $n_2 = 35000$ & $5 $ & $0.560$ & $7.05 \left(10 \right)^{-2}$ & $9.49 \left(10 \right)^{-2}$\\
  \hline
   $n_3 = 55000$ & $6$   & $0.548$ & $6.67 \left(10 \right)^{-2}$& $9.14 \left(10 \right)^{-2}$\\
  \hline
     $n_4 = 75000$ &$ 6$ & $0.532$& $6.24 \left(10 \right)^{-2}$ & $8.85 \left(10 \right)^{-2}$\\
  \hline
     $n_5 = 95000$  &$ 6$ & $0.512$ & $5.89 \left(10 \right)^{-2}$ & $8.47 \left(10 \right)^{-2}$\\
  \hline
     $n_6 = 115000$ &$ 6$ & $0.498$ & $5.62 \left(10 \right)^{-2}$ & $8.04 \left(10 \right)^{-2}$\\
  \hline
    $n_7 = 135000$ & $7$ & $0.495$ & $5.57 \left(10 \right)^{-2}$ & $7.66 \left(10 \right)^{-2}$\\
  \hline
      $n_8 = 155000$ &$ 7$ & $0.481$ & $5.28 \left(10 \right)^{-2}$ & $7.24 \left(10 \right)^{-2}$ \\
  \hline
    $n_9 = 175000$ & $7$ & $0.474$ & $5.01 \left(10 \right)^{-2}$ & $6.78 \left(10 \right)^{-2}$\\
  \hline
   $n_{10} = 195000$ & $7$  & $0.461$ & $4.90 \left(10 \right)^{-2}$ & $6.30 \left(10 \right)^{-2}$\\
  \hline
   $n_{11} = 215000$ & $7 $& $0.442$ & $4.69 \left(10 \right)^{-2}$ & $6.07 \left(10 \right)^{-2}$\\
  \hline
     $n_{12} = 235000$  &$ 7 $ & $0.425$ & $4.45 \left(10 \right)^{-2}$ & $5.82 \left(10 \right)^{-2}$\\
  \hline
     $n_{13} = 255000$ &$ 7$ & $0.411$ & $4.25 \left(10 \right)^{-2}$ & $5.54 \left(10 \right)^{-2}$\\
  \hline
    $n_{14} = 275000$ & $8$ & $0.408$ & $4.14 \left(10 \right)^{-2}$ & $5.16 \left(10 \right)^{-2}$ \\
  \hline
      $n_{15}= 295000$ & $8$  & $0.381$ & $4.09 \left(10 \right)^{-2}$ & $4.81 \left(10 \right)^{-2}$ \\
  \hline
   $n_{16} = 315000$ & $8$  & $0.360$ & $3.85 \left(10 \right)^{-2}$ & $4.53  \left(10 \right)^{-2}$\\
  \hline
   $n_{17} = 335000$ & $8 $ & $0.349$ & $3.56\left(10 \right)^{-2}$ & $4.29 \left(10 \right)^{-2}$\\
  \hline
     $n_{18} = 355000$ & $8 $& $0.330$ & $3.29\left(10 \right)^{-2}$  & $3.98 \left(10 \right)^{-2}$ \\
  \hline
     $n_{19} = 375000$  & $8$  & $0.320$   & $2.90 \left(10 \right)^{-2}$ & $3.75 \left(10 \right)^{-2}$ \\
  \hline
     $n_{20} = 395000$ &$ 8$ & $0.318$ & $2.62 \left(10 \right)^{-2}$  & $3.44 \left(10 \right)^{-2}$ \\
  \hline
\end{tabular}
\end{small}
  \label{A3:tab:Table11}
\end{table}

\textcolor{Crimson}{\section*{\textbf{Acknowledgments}}}

\textcolor{Aquamarine}{\textbf{This work has been supported in part by project MTM2015--71839--P (co-funded by Feder funds), of the DGI, MINECO, Spain.}}

\vspace{0.5cm}
\renewcommand\bibname{\textcolor{Crimson}{\textit{\textbf{References}}}}

\bibliographystyle{dinat}
\bibliography{Biblio}

\end{document}